\documentclass[a4paper]{article}
\usepackage{amssymb}
\usepackage{amsthm}
\usepackage{latexsym}
\usepackage{amsmath}
\usepackage{color}

\usepackage{comment}
\excludecomment{note}

\usepackage[margin=1in]{geometry}	
\usepackage{setspace}			

\newif\ifcol
\coltrue
\ifcol
\newcommand{\colorr}{\color[rgb]{0.8,0,0}}
\newcommand{\colorg}{\color[rgb]{0,0.5,0}}
\else
\newcommand{\colorr}{\color{black}}
\newcommand{\colorg}{\color{black}}
\fi

\newtheorem{definition}{Definition}[section]
\newtheorem{lemma}[definition]{Lemma}
\newtheorem{proposition}[definition]{Proposition}
\newtheorem{theorem}[definition]{Theorem}

\newtheorem{corollary}[definition]{Corollary}

\makeatletter
\@addtoreset{equation}{section}
\makeatother

\newcommand{\mca}{\mathcal{A}}
\newcommand{\mcb}{\mathcal{B}}
\newcommand{\mcc}{\mathcal{C}}
\newcommand{\mcd}{\mathcal{D}}
\newcommand{\mce}{\mathcal{E}}
\newcommand{\mcf}{\mathcal{F}}
\newcommand{\mcg}{\mathcal{G}}
\newcommand{\mci}{\mathcal{I}}
\newcommand{\mcl}{\mathcal{L}}

\newcommand{\mcn}{\mathcal{N}}

\newcommand{\mcr}{\mathcal{R}}

\newcommand{\mcu}{\mathcal{U}}
\newcommand{\mcv}{\mathcal{V}}
\newcommand{\mcx}{\mathcal{X}}

\newcommand{\mbbg}{\mathbb{G}}
\newcommand{\mbbh}{\mathbb{H}}
\newcommand{\mbbn}{\mathbb{N}}

\newcommand{\mbbr}{\mathbb{R}}

\newcommand{\mbbz}{\mathbb{Z}}
\newcommand{\mfa}{\mathfrak{A}}
\newcommand{\mfi}{\mathfrak{I}}
\newcommand{\mfs}{\mathfrak{S}}
\newcommand{\nn}{\nonumber}

\newcommand{\EQ}[1]{\begin{equation}{#1}\end{equation}}
\newcommand{\EQQ}[1]{\begin{equation}{#1 \nonumber}\end{equation}}
\newcommand{\EQN}[1]{\begin{equation}\begin{split}{#1}\end{split}\end{equation}}
\newcommand{\EQNN}[1]{\begin{equation}\begin{split}{#1 \nonumber}\end{split}\end{equation}}
\newcommand{\MAT}[2]{\left(\begin{array}{#1} #2 \end{array} \right)} 
\newcommand{\DIV}[2]{\left\{\begin{array}{#1} #2 \end{array} \right.} 

\newcommand{\PA}{\partial_\alpha}
\newcommand{\PB}{\partial_\gamma}

\newcommand{\TOD}{\overset{d}\to}
\newcommand{\TOP}{\overset{P}\to}

\newcommand{\EK}{\mce_{(k)}}
\newcommand{\EKK}{\mce_{(k')}}

\newcommand{\HAN}{\hat{\alpha}_n}
\newcommand{\HBN}{\hat{\gamma}_n}

\newcommand{\tmcd}{\tilde{\mcd}}
\newcommand{\FCSK}{\mcg_{k-1}}

\newcommand{\SUMM}{\sum_{m=1}^2}
\newcommand{\SUMMM}{\sum_{m_1,m_2=1}^2}
\newcommand{\SUMLL}{\sum_{l_1,l_2=1}^2}
\newcommand{\SUMP}{\sum_{p=0}^\infty}
\newcommand{\SUMPQ}{\sum_{p,q=0}^\infty}
\newcommand{\SUMIIp}{\sum_{\substack{i\leq {\bf M}_{1,n} \\ i'>{\bf M}_{1,n}}}}
\newcommand{\M}{{\bf M}}
\newcommand{\AS}{\quad {\rm a.s.}}
\newcommand{\TKI}{[t_{k-1}^n,t_k^n)}
\newcommand{\CFK}{\check{F}_{l_1,l_2,m_1,m_2, p,q,n,k}}
\newcommand{\DFK}{\dot{F}_{l_1,l_2,m_1,m_2,p,q,n,k}}
\newcommand{\DDFK}{\ddot{F}_{l_1,l_2,p,n,k}}
\newcommand{\CFN}{\check{F}_{l_1,l_2,m_1,m_2, p,q,n}}
\newcommand{\DFN}{\dot{F}_{l_1,l_2,m_1,m_2,p,q,n}}
\newcommand{\DDFN}{\ddot{F}_{l_1,l_2,p,n}}

\newcommand{\TAK}{\tilde{A}_{k,\epsilon}^{p,q}}

\newcommand{\al}{\alpha}
\newcommand{\del}{\delta}
\newcommand{\ep}{\epsilon}

\newcommand{\D}{\Delta}
\newcommand{\p}{\partial}

\newcommand{\Sig}{\Sigma}

\newcommand{\gam}{\gamma}

\newcommand{\aes}{\hat{\alpha}_{n}}

\newcommand{\ges}{\hat{\gamma}_{n}}

\newcommand{\diag}{{\rm diag}}

\title{Gaussian Quasi-Likelihood Estimation for Nonsynchronously Observed L\'evy processes}
\date{\today}
\author{Teppei Ogihara$^*$ and Yuma Uehara$^\dagger$\\
\begin{small}$*$ The University of Tokyo \quad $^\dagger$ Kansai University \end{small}
}

\begin{document}

\maketitle

\noindent
{\bf Abstract.}
We study parametric inference for a two-dimensional simple L\'evy process observed at high frequency over an expanding time horizon, where the two components are sampled in a possibly non-synchronous manner.
Combining Gaussian quasi-likelihood estimation for L\'evy SDEs with the non-synchronous quasi-likelihood methodology, we construct a Gaussian-type quasi-likelihood based on the covariance structure of non-synchronous increments, and we prove consistency and asymptotic normality of the resulting estimators.
The convergence rates coincide with those in the synchronous setting, while the asymptotic variance depends essentially on the sampling design and becomes more complicated under non-synchronicity.
To handle the different moment structure of pure-jump L\'evy processes, we introduce an additional convergence assumption and provide tractable sufficient conditions linked to the observation scheme.
Numerical experiments illustrate the finite-sample performance of the proposed estimators.

~


\section{Introduction}

Let $(\Omega, \mcf, P)$ be a probability space with a right-continuous filtration ${\bf F}=\{\mcf_t\}_{t\in \mbbr_+}$, where $\mbbr_+$ is the set of all nonnegative real numbers.
Let $X=(X^1_t,X^2_t)_{t\in\mbbr_+}$ be a $2$-dimensional stochastic process that is a solution to the stochastic differential equation
\EQ{dX_t=A(\alpha_0)dt+C(\gamma_0)dJ_t,}
where $J$ is a $2$-dimensional centered pure-jump L\'evy process whose L\'evy measure is denoted by $\nu$, and for bounded convex domains $\Theta_\alpha\subset \mbbr^{p_\alpha}$ and $\Theta_\gamma\subset \mbbr^{p_\gamma}$, 
$\alpha_0\in \Theta_\alpha$, $\gamma_0\in\Theta_\gamma$, and the coefficients $A:\Theta_\alpha\to\mbbr^2$, $C: \Theta_\gamma\to \mbbr^2\otimes \mbbr^2$ are measurable functions.
We hereafter write $\theta=(\alpha,\gamma)$ and $\Theta=\Theta_\alpha\times \Theta_\gamma$.
In this paper, we consider the estimation problem of $\theta$ under the situation where the elements $X^1$ and $X^2$ are observed at high frequency over an extended period, but their observation times may be different, that is, they are observed in a nonsynchronous manner.

In practice, we often need to take a nonsynchronous observation problem into account; typically, actual transaction data are irregularly recorded which makes it difficult to analyze the relevance of different stocks in the same way as in synchronous manner.
Such a problem has been discussed for a long time, and especially for diffusion processes, Hayashi and Yoshida (\cite{hay-yos05}, \cite{hay-yos08}, \cite{hay-yos11}) and Malliavin and Mancino (\cite{mal-man02}, \cite{mal-man09}) independently constructed consistent estimators of the  quadratic covariation.
These results were extended to work under the existence of market microstructure noise (\cite{ait-etal10, bar-etal11, bib-etal14, chr-etal10}).
Additionally, with the presence of jumps, \cite{Koi16} also constructed consistent estimators of the continuous-part integrated covariance based on threshold technique.
In the parametric estimation context, \cite{ogi-yos14} constructed a quasi maximum likelihood estimator and a Bayes-type estimator and showed its asymptotic mixed normality and tail probability estimates for nonsynchronously observed diffusion processes.
The aforementioned results consider an infill asymptotics over a fixed time interval.
As for the case where the sampling period is extending, \cite{ogi23} also clarified the asymptotic behavior of quasi maximum likelihood estimator.
Parametric estimation for infinite activity L\'evy processes and L\'evy driven stochastic differential equations has also been extensively studied under high-frequency sampling schemes.
For parametric L\'evy processes, Li et al.~\cite{li-etal24} studied Student L\'evy processes, Kawai and Masuda~\cite{kaw-mas13} investigated normal inverse Gaussian processes, and Masuda~\cite{mas09} and Brouste and Masuda~\cite{bro-mas18} treated stable processes.
For L\'evy driven stochastic differential equations, Cl\'ement and Gloter~\cite{cle-glo15,cle-glo19,cle-glo20} and Masuda~\cite{mas19} developed asymptotic theories for parametric estimation.
These works provide important foundations for high-frequency inference for jump-type models.
However, to the best of the authors' knowledge, there is no literature which addresses the parametric estimation problem of L\'{e}vy driven model with nonsynchronous data, and this motivates our study.

In this paper, we construct the quasi-log-likelihood function in the same manner as in \cite{ogi23}.
In the synchronous case, such a local Gaussian approximation enables us to estimate the drift and scale parameter with some loss of asymptotic efficiency, for more details, see \cite{mas13} and \cite{mas-ueh17}.
We show that our estimator still has the consistency and asymptotic normality as well as the synchronous case.
The convergence rate is also the same on one hand and the asymptotic variance differs and takes a more complex form essentially depending on the sampling design on the other hand.
Since the moment structure of pure jump L\'{e}vy process is quite different from that of Wiener process, we need to impose an additional convergent assumption (Assumption (A6)) for deducing our asymptotic results.
Although this assumption is closely related to the sampling condition and seems hard to check, we also study their tractable sufficient conditions.

In Masuda~\cite{mas13}, the estimator is constructed by approximating the covariance of the increments $\Delta X$ and using the likelihood function of a normal distribution with the same covariance. On the other hand, the quasi-likelihood functions of nonsynchronous observations of diffusion processes are constructed in Ogihara and Yoshida~\cite{ogi-yos14} and Ogihara~\cite{ogi23}, approximating the covariance of the increments $\Delta X$ from nonsynchronous observations and utilizing Gaussian likelihood functions. Fusing these two approaches, we also construct a Gaussian-type quasi-likelihood for nonsynchronous observations of the L\'evy process by calculating the covariance of $\Delta X$.

The remainder of the paper is organized as follows.
Section~\ref{main-section} introduces the observation scheme, constructs the Gaussian-type quasi-likelihood and the estimators, and states the assumptions and the main asymptotic normality theorem.
Section~\ref{suff-cond-section} provides sufficient conditions for the sampling assumptions.
Section~\ref{num-section} presents a blocking implementation of the quasi-likelihood and reports numerical experiments.
Section~\ref{proofs-section} contains the proofs of the main results.


\section{Main results}\label{main-section}

Let $(h_n)_{n=1}^\infty$ be a sequence of positive numbers satisfying
\begin{equation}\label{bcond}
    h_n\to 0, \quad n^{1-\frac{1}{m}}h_n\to \infty,\quad nh_n^2\to 0,
\end{equation}
as $n\to\infty$ for some integer $m\geq3$.
Intuitively, $n$ is of the order of the number of observations and $h_n$ is of the order of the length of the observation intervals.
More precise assumptions for the observations are given later.
Let $\mathbb{N}$ and $\mathbb{Z}_+$ denote the sets of all positive integers and all nonnegative integers, respectively.
Let $T_n=nh_n$.
For a $k\times l$ matrix $M$ and a vector $v$ of size $k$, 
$[M]_{ij}$ denotes $(i,j)$ element of $M$, $\lVert M\rVert$ denotes the operator norm induced by the Euclidean norm, and $[v]_i$ denotes $i$-th element of $v$.
For a vector $x=(x_1,\cdots, x_k)$ and positive integer $l$, let $\partial_x^l=(\frac{\partial^l}{\partial x_{i_1}\cdots \partial x_{i_l}})_{i_1,\cdots,i_l=1}^k$. $\partial_x^0$ denotes the identity.

For $l\in\{1,2\}$ and $n\in\mathbb{N}$, ${\bf M}_{l,n}$ is a random positive integer, $\{S^{n,l}_i\}_{i=0}^{{\bf M}_{l,n}}$ are strictly increasing random times satisfying $S_0^{n,l}=0$ and $S_{\M_{l,n}}^{n,l}\leq T_n$.
We assume that $\{S_i^{n,l}\}_{i,l}$ are independent of $(J_t)_{t\in \mbbr^+}$.
We will assume in (A3) that $T_n-S_{{\bf M}_{l,n}}^{n,l}\TOP 0$ for $l\in\{1,2\}$.

Let ${\bf M}_n={\bf M}_{1,n}+{\bf M}_{2,n}$, and let $I^l_i=[S^{n,l}_{i-1},S^{n,l}_i)$. For $1\leq i\leq {\bf M}_n$, let 
\EQQ{\varphi(i)=\DIV{ll}{i, & {\rm if} \ i\leq {\bf M}_{1,n},\\ i-{\bf M}_{1,n}, & {\rm if} \ i>{\bf M}_{1,n},}
\quad \psi(i)=\DIV{ll}{1, & {\rm if} \ i\leq {\bf M}_{1,n},\\ 2, & {\rm if} \ i>{\bf M}_{1,n}.}}
We set $I_i=I^{\psi(i)}_{\varphi(i)}$ for $1\leq i\leq {\bf M}_n$.
We often consider a vector of the observations of $X^1$ and $X^2$ in a row, and therefore, the following notation is convenient.
For a two-dimensional stochastic process $(U_t)_{t\geq 0}=((U_t^1,U_t^2)^\top)_{t\geq 0}$,
let $\Delta_{i,t}^l U=U^l_{t\wedge S^{n,l}_i}-U^l_{t\wedge S^{n,l}_{i-1}}$, and let $\Delta^l_t U=(\Delta_{i,t}^l U)_{1\leq i\leq {\bf M}_{l,n}}$.
Let $\Delta_t U=((\Delta^1_t U)^\top, (\Delta^2_t U)^\top)^\top$.
Let $\Delta_i^l U=\Delta_{i,T_n}^l U$, and $\Delta U=\Delta_{T_n} U$.

Let $V^l(\gamma)=[CC^\top(\gamma)]_{ll}$, $V^{1,2}(\gamma)=[CC^\top(\gamma)]_{12}$,
and $\mcv_t^l(\alpha)=tA^l(\alpha)$.
\begin{note}
    {\colorr Change notation if $C$ will be not depending on $X$}
\end{note}
Let
\EQQ{S_n(\gamma)=\MAT{cc}{{\rm diag}((V^1(\gamma)|I_i^1|)_i) & (V^{1,2}(\gamma)|I_i^1\cap I_j^2|)_{ij} \\ 
(V^{1,2}(\gamma)|I_i^1\cap I_j^2|)_{ij}^\top & {\rm diag}((V^2(\gamma)|I_j^2|)_j) }.}
\begin{note}
{\colorr Explain the reason of th invertivity of $S_n$}
\end{note}
We define a quasi-log-likelihood function for $\gamma$ by
\EQQ{\mbbh_n^1(\gam)=-\frac{1}{2}\D X^\top S_n^{-1}(\gam)\D X-\frac{1}{2}\log\det S_n(\gam)}
Then, a maximum-likelihood-type estimator is defined by
\EQQ{\HBN\in {\rm argmax}_{\gam\in \bar{\Theta}_\gam}\mbbh_n^1(\gam).}

Unlike in the case of diffusion processes, to estimate drift parameter, we consider the one-step estimation procedure for some technical reasons.
As an initial estimator, we adopt the following least square estimator:
\begin{align*}
\tilde{\alpha}_n\in\text{argmin}_{\alpha\in \bar{\Theta}_\al} \mbbg_{2,n}(\alpha),
\end{align*}
where
\begin{align*}
    \mbbg_{2,n}(\alpha)&=(\D X-\D \mcv(\al))^\top\diag\{(|I_i^1|)_i,(|I_j^2|)_j\}^{-1}(\D X-\D \mcv(\al))\\
    &=\sum_{i\leq {\bf M}_{1,n}} \frac{[\D X-\D \mcv(\al)]_i^2}{|I_i^1|}+\sum_{j> {\bf M}_{1,n}} \frac{[\D X-\D \mcv(\al)]_j^2}{|I_{j-\M_{1,n}}^2|}.
\end{align*}
Then, we define our one-step estimator $\aes$ as 
$$\aes=\tilde{\al}_n-\Big(\p_\al^2\mbbh_n^2(\tilde{\al}_n)\Big)^{-1}\p_\al \mbbh^2_n(\tilde{\al}_n),$$
where the quasi-likelihood function $\mbbh_n^2$ is defined by
$$\mbbh_n^2(\al)=-\frac{1}{2}(\D X-\D \mcv(\al))^\top S_n^{-1}(\HBN)(\D X-\D \mcv(\al)).$$

The existence of the inverse of $S_n$ is ensured under Assumption (A2) in the next subsection by Proposition 1 in Section 2 of Ogihara and Yoshida (2014).


For a matrix $A$, $A^\top$ denotes the transpose, and we denote $A^{\otimes 2}=AA^\top$.
Let $I$ denote the unit matrix of suitable size.

\begin{description}
    \item[Assumption (A1).]
The process $J$ satisfies 
$$E[J_1]=0, \ E[J_1^{\otimes 2}]=I, \ E[|J_1|^q]<\infty,$$
for all $q>0$.
\end{description}
From L\'{e}vy-Khintchine representation (cf. \cite{App09}), we can easily deduce that $E[J_t]=0$ and $E[J_t^{\otimes 2}]=tI$ for any $t>0$.

~

For the integer $m$ introduced in \eqref{bcond}, let $\tilde{m}=m\vee 4$.
\begin{description}
\item[Assumption (A2).] 
For all nonnegative integer $k\leq \tilde{m}$, $\PA^kA(\alpha)$ and $\PB^kC(\gamma)$ exist and are continuous on $\bar{\Theta}_\alpha$ and $\bar{\Theta}_\gamma$, respectively.
There exists a positive constant $c_1$ such that $c_1I\leq CC^\top(\gamma)$ for any $\gamma\in \bar{\Theta}_\gamma$.
\end{description}

Let 
\EQNN{r_n&= \max_{l\in\{1,2\}}\max_{1\leq i\leq{\bf M}_{l,n}}\left\{(S_i^{n,l}-S_{i-1}^{n,l})\;\vee\;(T_n-S_{{\bf M}_{l,n}}^{n,l})\right\}, \\
\underbar{r}_n &= \min_{l\in\{1,2\}}\min_{1\leq i\leq{\bf M}_{l,n}}(S_i^{n,l}-S_{i-1}^{n,l}), }
and let
\EQQ{G=\bigg\{\frac{|I^1_i\cap I^2_j|}{|I_i^1|^{1/2}|I_j^2|^{1/2}}\bigg\}_{1\leq i\leq {\bf M}_{1,n}, 1\leq j\leq {\bf M}_{2,n}}.}

\begin{description}
\item[Assumption (A3).] 
$r_n\TOP 0$ as $n\to\infty$, and there exists $\delta >0$ such that $n^\delta T_n^{-1}(h_n\underbar{r}_n^{-1})^4\TOP 0$.
\end{description}

Let $\mfs$ be the set of all partitions $(s_k)_{k=0}^\infty$ of $[0,\infty)$ satisfying $\sup_{k\geq 1}|s_k-s_{k-1}|\leq 1$ and $\inf_{k\geq 1}|s_k-s_{k-1}|>0$.
For $(s_k)_{k=0}^\infty\in \mfs$, let $M_{l,k}=\#\{i;\sup I_i^l\in (s_{k-1},s_k]\}$ and $q_n=\max\{k;s_k< nh_n\}$,
and let $\EK^l$ be an $\M_{l,n}\times \M_{l,n}$ matrix satisfying
$[\EK^l]_{ij}=1$ if $i=j$ and $\sup I_i^l\in (s_{k-1},s_k]$, and otherwise $[\EK^l]_{ij}=0$.

\begin{description}
\item[Assumption (A4).] 
There exist positive constants $a_0^1$ and $a_0^2$ such that $\{h_nM_{l,q_n+1}\}_{n=1}^\infty$ is $P$-tight and
\EQQ{\max_{1\leq k\leq q_n} |h_nM_{l,k}-a_0^l(s_k-s_{k-1})|\TOP 0}
as $n\to\infty$ for $l\in \{1,2\}$ and any partition $(s_k)_{k=0}^\infty \in\mfs$. 
Moreover, for any $p\in \mbbn$, there exists a nonnegative constant $a_p^1$ such that 
\EQQ{\max_{1\leq k\leq q_n}|h_n{\rm tr}(\EK^1(GG^\top)^p)-a_p^1(s_k-s_{k-1})|\TOP 0}
as $n\to\infty$ for any partition $(s_k)_{k=0}^\infty \in\mfs$.
\end{description}
In particular, by setting $s_k=k$ for $k\geq 0$, $q_n$ is the largest integer strictly smaller than $nh_n$, and Assumption (A4) implies that 
$$\frac{1}{nh_n}\bigg|\sum_{l=1}^2\sum_{k=1}^{q_n}(h_nM_{l,k}-a_0^l(s_k-s_{k-1}))\bigg| \TOP 0,$$
and consequently we obtain
\EQ{\label{Mn-order-eq} \frac{{\bf M}_n}{n}=\frac{1}{nh_n}\sum_{l=1}^2\sum_{k=1}^{q_n+1}h_nM_{l,k}\TOP a_0^1+a_0^2.}

Let $\mfi_l=(|I_i^l|^{1/2})_{i=1}^{\M_{l,n}}$.
\begin{description}
\item[Assumption (A5).] For $p\in \mbbz_+$, there exist nonnegative constants $f_p^{1,1}$, $f_p^{1,2}$, and $f_p^{2,2}$ such that 
$\{|\mce_{(q_n+1)}^l\mfi_l|\}_{n=1}^\infty$ is $P$-tight for $l\in\{1,2\}$, and
\EQNN{\max_{1\leq k\leq q_n}|\mfi_1^\top \EK^1(GG^\top)^p\mfi_1-f_p^{1,1}(s_k-s_{k-1})|&\TOP 0, \\
\max_{1\leq k\leq q_n}|\mfi_1^\top \EK^1(GG^\top)^pG\mfi_2-f_p^{1,2}(s_k-s_{k-1})|&\TOP 0, \\
\max_{1\leq k\leq q_n}|\mfi_2^\top \EK^2(G^\top G)^p\mfi_2-f_p^{2,2}(s_k-s_{k-1})|&\TOP 0
}
as $n\to\infty$ for any partition $(s_k)_{k=0}^\infty \in\mfs$.
\end{description}

Assumption (A4) corresponds to [A3$'$] in Ogihara and Yoshida~\cite{ogi-yos14}.
Though it is difficult to directly check (A4) and (A5) for general sampling schemes, 
tractable sufficient conditions for these conditions are studied in Section~2.4 in~\cite{ogi23}.

Let 
$$\chi_{Cz}(t)=([C(\gamma_0)z]_1(1_{I_i^1}(t))^\top_{1\leq i\leq {\bf M}_{1,n}}, [C(\gamma_0)z]_2(1_{I_i^2}(t))^\top_{1\leq i\leq {\bf M}_{2,n}})^\top,$$
and we define
$$\mathcal{S}(t,z)=\MAT{c}{(h_n/2\chi_{Cz}(t)^\top \partial_{\gamma_i}S_n^{-1}(\gamma_0) \chi_{Cz}(t))_{i=1}^{p_\gamma} \\ 
(- \partial_{\alpha_j}\Delta \mcv(\alpha_0)^\top S_n^{-1}(\gamma_0) \chi_{Cz}(t))_{j=1}^{p_\alpha}},$$
and let
$$F_n =\frac{1}{T_n} \int_0^{T_n}\int \mathcal{S}\mathcal{S}^\top(t,z) \nu(dz)dt.$$

\begin{description}
\item[Assumption (A6).]\label{convav}
There exists a matrix $\Gamma_1$ 
such that $F_n \TOP \Gamma_1$ as $n\to\infty$.
\end{description}
Since $F_n$ is nonnegative definite, $\Gamma_1$ is nonnegative definite.
\begin{note}
    {\colorr explanation: A6 is a new condition, sufficient conditions for (A6) is given in Sec X (the same as sufficient conditions in Section 2.4 of~\cite{ogi23}), etc.
    
$$\dot{F}_{1,1,1,1,p,q,n}=\frac{h_n}{n}\sum_{i,i'}|I_i\cap I_{i'}|\frac{[(GG^\top)^p]_{ii}[(GG^\top)^q]_{ii}}{|I_i|^2}
=\frac{h_n}{n}\sum_i\frac{[(GG^\top)^p]_{ii}[(GG^\top)^q]_{ii}}{|I_i|}.$$
For synchronous observations, since $G=I$,
$$\dot{F}_{1,1,1,1,p,q,n}=\frac{h_n}{n}\sum_i|I_i|^{-1}.$$
}
\end{note}

\begin{description}
\item[Assumption (A7).] 
The constant $a_1^1$ in (A4) is positive, and there exist positive constants $c_2$ and $c_3$ such that 
\EQQ{\lVert CC^\top(\gamma) -CC^\top(\gamma_0)\rVert^2\geq c_2|\gamma-\gamma_0|^2, 
\quad |A(\alpha)-A(\alpha_0)|^2\geq c_3|\alpha-\alpha_0|^2
}
for any $\gamma\in \Theta_\gamma$ and $\alpha\in \Theta_\alpha$.
\end{description}
Assumption (A7) will be used to identify the parameters $\gamma$ and $\alpha$ from the data.

Following the information-matrix construction in Ogihara~\cite{ogi23},
we define the matrix $\Gamma_2$ for the present constant-coefficient model.
For $l\in\{1,2\}$, let
\[ \rho(\gamma)=\frac{V^{1,2}(\gamma)}{\sqrt{V^1(\gamma)V^2(\gamma)}},
 \quad \rho_0=\rho(\gamma_0),
 \quad B_l(\gamma)=\sqrt{\frac{V^l(\gamma_0)}{V^l(\gamma)}},
 \quad 
 \phi_l(\alpha)=\frac{A^l(\alpha)-A^l(\alpha_0)}{\sqrt{V^l(\gamma_0)}}.
 \]
For a column vector $x$, write
$x^{\otimes2}=xx^\top$. 
Using the constants in (A4) and (A5), set
$a(\rho)=\sum_{p=1}^{\infty}a_p^1\rho^{2p}$ for $|\rho|<1$.
Let
\EQNN{\mathcal I_\gamma &= 2\sum_{l=1}^2\bigl(a_0^l+a(\rho_0)\bigr) (\partial_\gamma B_l)^{\otimes2}(\gamma_0)
 +\frac{a'(\rho_0)}{\rho_0}(\partial_\gamma \rho)^{\otimes2}(\gamma_0)
 -a(\rho_0)
\left(\frac{\partial_\gamma \rho}{\rho_0} - \partial_\gamma B_1 - \partial_\gamma B_2 \right)^{\otimes2}(\gamma_0), \\
 \mathcal I_\alpha &= \sum_{p=0}^{\infty} \rho_0^{2p}
     \left\{\sum_{l=1}^2 f_p^{l,l}(\partial_\alpha \phi_l)^{\otimes2}(\alpha_0)
       -\rho_0 f_p^{1,2}
\bigl(\partial_\alpha \phi_1\partial_\alpha \phi_2^\top+\partial_\alpha \phi_2\partial_\alpha \phi_1^\top\bigr)(\alpha_0) \right\}.
}
At $\rho_0=0$, this definition is understood by continuous extension.
Finally, let
\[
 \Gamma_2
   =\begin{pmatrix}
        \mathcal I_\gamma&0\\
        0&\mathcal I_\alpha
      \end{pmatrix}.
\]
Under (A1)--(A7), $\Gamma_2$ is symmetric and
positive definite, by the information bounds in Ogihara~\cite{ogi23} (see (3.39) and the proof of Theorem 2.3)
and the identifiability condition (A7).

\begin{theorem}\label{asym-norm-thm}
Assume (A1)--(A7). Then, 
$$
\sqrt{T_{n}}\binom{\hat{\gamma}_{n}-\gamma_{0}}{\hat{\alpha}_{n}-\alpha_{0}} \TOD N\left(0, \Gamma_{2}^{-1} \Gamma_{1} \Gamma_{2}^{-1}\right)
$$
as $n \rightarrow \infty$.

\end{theorem}

This result shows that nonsynchronicity does not reduce the
$\sqrt{T_n}$ convergence rate
for synchronous Gaussian quasi-likelihood estimation of L\'evy-driven
models in Masuda~\cite{mas13}. The lengths and overlaps of the observation
intervals enter the limiting covariance matrix. In the sandwich covariance
$\Gamma_2^{-1}\Gamma_1\Gamma_2^{-1}$, $\Gamma_2$ represents the limiting
curvature of the appropriately normalized quasi-log-likelihoods, whereas
$\Gamma_1$ is the limiting covariance of the corresponding normalized
score vector. The latter also depends on the jump distribution through
the L\'evy measure.

Since the Gaussian criterion is a quasi-likelihood, the usual information
identity need not hold, and $\Gamma_1$ and $\Gamma_2$ generally differ.
Neither the common convergence rate nor the sandwich representation
implies asymptotic efficiency or equality of precision with synchronous
sampling. Thus, for the sampling schemes covered by the theorem,
irregular and nonsynchronous observations remain compatible with
estimation at the $\sqrt{T_n}$ rate, while the construction of asymptotic
standard errors must account for the sampling design and the jump
distribution.

\section{Sufficient conditions for the assumptions}\label{suff-cond-section}

It is not easy to check whether a random sampling scheme satisfies Assumptions (A4)--(A6). In this section, we give sufficient conditions for (A4)--(A6) using an approach similar to that of \cite{ogi23}.

Let $\mcn_t^{n,l}=\sum_{i=1}^{{\bf M}_{l,n}}1_{\{S^{n,l}_i\leq t\}}$. For any $q>0$, we consider the following conditions.

\begin{description}
\item[Assumption (B1-$q$).]
$$\sup_{n\geq 1}\max_{l\in\{1,2\}}\sup_{0\leq t\leq (n-1)h_n}E[(\mcn_{t+h_n}^{n,l}-\mcn_t^{n,l})^q]<\infty.$$
\item[Assumption (B2-$q$).]
$$\limsup_{u\to\infty}\sup_{n\geq 1}\max_{l\in\{1,2\}}\sup_{0\leq t\leq nh_n-uh_n}u^qP(\mcn_{t+uh_n}^{n,l}-\mcn_t^{n,l}=0)<\infty.$$
\item[Assumption (B3).] There exist $\epsilon\in (0,1)$ and partitions $(t_k^n)_{k=0}^{L_n}$ of $[0,T_n]$ for each $n\in \mbbn$ ($0=t_0^n<t_1^n<\cdots < t_{L_n}^n=T_n$) such that $L_n\to \infty$, $n^{\epsilon-1}L_n\to 0$, $L_n(nh_n)^{-1}E[r_n^4]^{1/4}\to 0$ as $n\to \infty$, 
\EQ{\label{B3-eq} \sup_n\frac{\max_{1\leq k\leq L_n}(t_k^n-t_{k-1}^n)}{\min_{1\leq k\leq L_n}(t_k^n-t_{k-1}^n)}<\infty, \quad \sup_{n,k} E\bigg[\bigg(\frac{h_n^2}{(t_k^n-t_{k-1}^n)}\sum_{\substack{1\leq i\leq \M_n \\ I_i\cap 
 \TKI\neq \emptyset}}|I_i|^{-1}\bigg)^4\bigg]<\infty.}
\end{description}

Assumption (B3) is not necessary in the previous work Ogihara~\cite{ogi23} for nonsynchronously observed diffusion processes.
Therefore, this condition is specific to the nonsynchronous observation of L\'evy processes.
A well-known example for $\mcn_t^{n,l}$ is the Poisson process, in which the condition $E[|I_i^l|^{-1}]=\infty$ holds for each $i,l$, and thus (B3) is not satisfied.

We examine why this condition is necessary.
To establish asymptotic normality of the estimator, one needs the limit
  \[ \lim_{n\to\infty} E\Bigl[\bigl(\sqrt{h_n/n}\,\PB \mbbh_n^1 (\gam_0) \bigr)^2\Bigr]
  \]
  to be finite, as discussed in Section~\ref{asym-norm-subsection}.

For simplicity, consider the one-dimensional case with the drift $\mathcal V\equiv 0$, so that
  \[
    S_n = \diag(\bigl(V^1\,|I_i^1|\bigr)_i).
  \]
Then we have
\EQQ{\PB \mbbh_n^1 = - \frac{\PB (V^1)^{-1}}{2} \sum_{i=1}^{\M_{1,n}}\frac{(\Delta_i^1 X)^2}{|I_i^1|}-\frac{\M_{1,n}}{2} \frac{\PB V^1}{V^1},}
and hence, we obtain
  \[
    \frac{h_n}{n}E\bigl[(\PB \mbbh_n^1)^2\mid \Pi \bigr]
      = \frac{h_n}{4n} \bigl(\PB ((V^1)^{-1})\bigr)^2
         E\Bigl[\Bigl(\sum_{i=1}^{\M_{1,n}} \frac{(\Delta_i^1 X)^2}{|I_i^1|}\Bigr)^2 \Big| \Pi \Bigr]
        - \frac{h_n}{4n}\Bigl(\frac{\PB V^1}{V^1}\Bigr)^2 \M_{1,n}^2.
  \]
  The expectation of the first term on the right--hand side is
\EQN{\label{invI-eq} E\Bigl[\Bigl(\sum_{i=1}^{\M_{1,n}} \frac{(\Delta_i^1 X)^2}{|I_i^1|}\Bigr)^2 \Big| \Pi \Bigr]
=E\bigg[\sum_{i=1}^{\M_{1,n}} \frac{(\Delta_i^1 X)^4}{|I_i^1|^2} \Big| \Pi \bigg]+E\bigg[\sum_{i\neq j}\frac{(\Delta_i^1 X)^2(\Delta_j^1 X)^2}{|I_i^1||I_j^1|} \Big| \Pi \bigg].
}
For the L\'evy process $X$, the forth moment is calculated as 
  \[
    E\bigl[(\Delta_i^1 X)^4\mid \Pi \bigr] = c\,|I_i^1|+ ({\rm higher~order~terms}),
  \]
for some constant $c$, and hence
  \[
    E\left[\,\sum_i |I_i^1|^{-1}\right]
  \]
must be $O(n/h_n)$ so that the limit of the right-hand side of (\ref{invI-eq}) is finite
(In contrast, when $X$ is a Brownian motion,
\(E[(\Delta X_i^1)^4\mid \Pi] = c\,|I_i^1|^2+ ({\rm higher~order~terms})\), 
so (B3) is not needed).
Therefore, an assessment for the sum of $|I_i^1|^{-1}$ is required. Furthermore, the necessity of applying the local sum argument and evaluating the fourth moment leads to an assumption of the form (B3).


We consider mixing properties of $\mcn^{n,l}_t$. Let
$$\mathcal{H}_{i,j}^n=\sigma(\mcn_t^{n,l}-\mcn_s^{n,l};ih_n\leq s<t\leq jh_n, l\in \{1,2\})$$
for $0\leq i,j\leq n$, and let $\alpha_0^n=1/4$, and for $k\in \mbbn$,
$$\alpha_k^n=0\vee \sup_{\substack{1\leq i,j\leq n-1 \\ j-i\geq k}}\sup_{A\in \mathcal{H}_{0,i}^n, B\in \mathcal{H}_{j,n}^n}|P(A\cap B)-P(A)P(B)|.$$

\begin{description}
\item[Assumption (B4-$q$).]
\EQ{\label{alpha-k-eq} \sup_{n\geq \mbbn}\sum_{k=0}^\infty (k+1)^q\alpha_k^n<\infty.}
\item[Assumption (B5).]
There exists a matrix $\Gamma_1$ such that $E[F_n]\to \Gamma_1$
as $n\to\infty$.
\item[Assumption (B6).]
There exist positive constants $a_0^1$ and $a_0^2$, and nonnegative constants $a_{p+1}^1$, $f_p^{1,1}$, $f_p^{1,2}$, and $f_p^{2,2}$ for $p\in \mbbz_+$ such that
$\{E[h_nM_{l,q_n+1}]\}_{n=1}^\infty$ and $\{E[|\mce_{(q_n+1)}^l\mfi_l|]\}_{n=1}^\infty$ are bounded and
\EQN{\label{A4A5-average-cond} \max_{1\leq k\leq q_n}|h_nE[M_{l,k}]-a_0^l(s_k-s_{k-1})|&\to 0, \\
\max_{1\leq k\leq q_n}|h_nE[{\rm tr}(\EK^1(GG^\top)^p)]-a_p^1(s_k-s_{k-1})|&\to 0, \\
\max_{1\leq k\leq q_n}|E[\mfi_1^\top \EK^1(GG^\top)^p\mfi_1]-f_p^{1,1}(s_k-s_{k-1})|&\to 0, \\
\max_{1\leq k\leq q_n}|E[\mfi_1^\top \EK^1(GG^\top)^pG\mfi_2]-f_p^{1,2}(s_k-s_{k-1})|&\to 0, \\
\max_{1\leq k\leq q_n}|E[\mfi_2^\top \EK^2(G^\top G)^p\mfi_2]-f_p^{2,2}(s_k-s_{k-1})|&\to 0}
as $n\to\infty$ for $p\in\mbbz_+$, $l\in\{1,2\}$ and any partition $(s_k)_{k=0}^\infty\in\mfs$.
\end{description}
(B6) represents that (A4) and (A5) are reduced to the conditions for the expectations.

\begin{proposition}\label{A6-suff-prop}
Assume (A1), (A2), (B2-$r$), (B3), (B4-$1$) and (B5) hold for some $r>0$.
Then, (A6) holds.
\end{proposition}

Consequently, we obtain the following corollary from Proposition~\ref{A6-suff-prop} together with
Propositions 2.7--2.9, and 2.11 in Ogihara~\cite{ogi23}.
\begin{corollary}
Assume (A1), (A2), (B1-$q$), (B2-$q$), (B3), (B4-$q$), (B5) and (B6) for any $q>0$,
and there exists $\eta<1/(4m)$ such that 
\EQ{\label{min-r-eq} n^{-\eta} h_n\underbar{r}_n^{-1}\TOP 0 \quad (n\to\infty).}
Then, (A3)--(A6) hold. Moreover, $a_1^1$ in (A4) is positive.
\end{corollary}

~

Next, we provide a class of sampling schemes based on renewal processes for which (B1-$q$), (B2-$q$), (B3), (B4-$q$), (B5) and (B6) are satisfied.

\begin{description}
\item[Assumption (R).] 
For $l \in\{1,2\}$, let $\left(\zeta_{i}^{l}\right)_{i \in \mbbn}$ be independent, identically distributed random variables, satisfying $\underline{\zeta}\leq \zeta_i^l\leq \bar{\zeta}$ a.s.\ for all $i\in\mbbn$ with constants $0<\underline{\zeta}\leq \bar{\zeta}$.
The two sequences are mutually independent.
Assume that the distribution of $\zeta_1^l$ is not lattice.
Moreover, for each $l\in\{1,2\}$, there exist constants
$\underline{\zeta}\leq u_l<v_l\leq\bar{\zeta}$ and
$\eta_l\in(0,1)$ such that
\[
P(\zeta_1^l\in A)
\geq \frac{\eta_l}{v_l-u_l}
       \int_A 1_{[u_l,v_l]}(x)\,dx
\]
for every Borel set $A\subset\mathbb R$.

The observation times $\{S_i^{n,l}\}$ are given by 
\EQ{\label{renewal-def-eq} \bar{\mcn}_{t}^{l}=\sup \{k \in \mbbz_{+} ; \sum_{1 \leq i \leq k} \zeta_{i}^{l} \leq t\}, \quad \mcn_t^{n,l}=\bar{\mcn}_{h_n^{-1}t}^l,}
$$S_i^{n,l}=\inf\{t\geq 0; \mcn_t^{n,l}=i\}, \quad {\bf M}_{l,n}=\max\{i\geq0:S_i^{n,l}\leq T_n\}.$$
\end{description}

\begin{proposition}\label{renewal-exp-prop}
Assume (A1), (A2) and (R).
Then, (B1-$q$), (B2-$q$), (B3), (B4-$q$), (B5), (B6) hold for any $q>0$. As a consequence, (A3)--(A6) hold and $a_1^1$ in (A4) is positive.
\end{proposition}


\section{Numerical Experiment}\label{num-section}

Before introducing the statistical model for our simulation study, we explain the blocking technique which relieves a heavy computational load to calculate our quasi-likelihood function, especially the inverse matrix of $S_n$. 
The same technique used in \cite{ogi18} and was originally introduced to handle measurement error by \cite{GloJac01}:
Fix $b\in\mbbn$.
We divide the whole observation interval $[0,T_n]$ into $b$ equidistant intervals, and construct the quasi-likelihood functions over each interval; finally, the sum of these functions is used for our parameter estimation.
More specifically, for $l\in\{1,2\}$ and $k\in\{1,\dots,b\}$, let
\begin{align*}
    & B_k^{n,l}=\Bigg\{i: [S_{i-1}^{n,l},S_{i}^{n,l})\subset \Bigg[(k-1)\frac{T_n}{b},k\frac{T_n}{b}\Bigg)\Bigg\}, \ \Delta X_{B_k}=((\Delta_{i_1} X^1)_{i_1\in B_k^{n,1}}^\top,(\Delta_{i_2} X^2)_{i_2\in B_k^{n,2}}^\top)^\top,\\
    &\Delta \mcv_{B_k}=((\Delta_{i_1} \mcv^{1}(\al))_{i_1\in B_k^{n,1}}^\top,(\Delta_{i_2} \mcv^2(\al))_{i_2\in B_k^{n,2}}^\top)^\top,\\ 
    &S_{n,k}(\gam)=\MAT{cc}{{\rm diag}((V^1(\gamma)|I_{i}^1|)_i) & (V^{1,2}(\gamma)|I_{i}^1\cap I_{j}^2|)_{ij} \\ 
    (V^{1,2}(\gamma)|I_{i}^1\cap I_{j}^2|)_{ij}^\top & {\rm diag}((V^2(\gamma)|I_{j}^2|)_j) }, \quad i\in B_k^{n,1}, \ j\in B_k^{n,2}.
\end{align*}
Then, our estimators are defined by $\ges\in\text{argmax}_{\gam\in\Theta_\gam} \mbbh_{1,n}^B(\gam)$ and $\aes\in\text{argmax}_{\al\in\Theta_\al} \mbbh_{2,n}^B(\al)$, where
\begin{align*}
    &\mbbh_{1,n}^B(\gam)=-\sum_{k=1}^b\bigg[(\Delta X_{B_k})^\top S_{n,k}^{-1}(\gam)(\Delta X_{B_k})+\log \det S_{n,k}(\gam)\bigg],\\
    &\mbbh_{2,n}^B(\al)=-\sum_{k=1}^b(\Delta X_{B_k}-\Delta \mcv_{B_k}(\alpha))^\top S_{n,k}^{-1}(\ges)(\Delta X_{B_k}-\Delta \mcv_{B_k}(\alpha)).
\end{align*}

In the rest of this section, we consider the following two-dimensional stochastic process:
$$dX_t=\begin{pmatrix}\alpha_1\\\alpha_2\end{pmatrix}dt+\begin{pmatrix}\gamma_1&\gamma_3\\\gamma_3&\gamma_2\end{pmatrix}d\begin{pmatrix}J_{1,t}\\J_{2,t}\end{pmatrix},$$
where the one-dimensional L\'evy processes $J_{1,t}$ and $J_{2,t}$ are independent and both of them obey the normal inverse Gaussian distribution $NIG(3,0,3t,0)$ for each $t$.
The normal inverse Gaussian random variable is defined as the normal variance-mean mixture of the inverse Gaussian random variable, and the mean and variance of $NIG(\alpha', \beta, \delta,\mu)$ are given by $\mu+\delta\beta/\sqrt{(\alpha')^2-\beta^2}$ and $\frac{\del(\al')^2}{((\al')^2-(\beta)^2)^{3/2}}$, respectively (cf. \cite{Bar97}).
The parameter spaces are $\Theta_\al=(-2,2)\times(-2,2)$ and $\Theta_\gam=(0.7,3)\times(0.7,3)\times(-0.5,0)$, respectively.
We set the sampling times $\{S_i^{n,l}\}_{i=1}^{{\bf M}_{l,n}}$ as
$$S_i^{n,l}=\sum_{k=1}^i h_nY_k^{l}, \quad l=1,2$$
where $Y_1^l,Y_2^l,\dots,$ are i.i.d. random variables with $Y_1^l\sim \Gamma(4.5,3)$, and the random sequences $\{Y_k^1\}$ and $\{Y_k^2\}$ are supposed to be independent.
Hence, the number of samples over the interval $[0,t]$ becomes the following renewal process: 
$$\mcn_t^{n,l}=\sum_{i=1}^{{\bf M}_{l,n}}1_{\{S^{n,l}_i\leq t\}}, \quad l=1,2.$$
In such a setting, we independently generate 1000 paths of the solution process $X$ based in Euler-Maruyama scheme for the three cases: $(T_n,h_n,b)=(50,0.03,20), (100,0.01,40), (200,0.004,75)$.
In our simulation, the true values of each parameter are 
$$(\al_{1,0},\al_{2,0},\gam_{1,0},\gam_{2,0},\gam_{3,0})=(-0.4,0.5,2,1,-0.2),$$
and the initial values $\al_{\text{initial}}$ and $\gam_{\text{initial}}$ are set to $\al_{\text{initial}}=(-1,1)$ and $\gam_{\text{initial}}=(2,1,-0.2)$, respectively.
For the calculation of the estimators, we use the blocking technique introduced at the beginning of this section.
Our simulation result is given in Table \ref{tab1}.
The numerical result seems to reflect the theoretical result given in this paper.


\begin{table}\label{sres}
    \centering
    \begin{tabular}{cccccccc}\hline
         $T_n$& $h_n$ & $b$ & $\hat{\gam}_{1,n}$ &$\hat{\gam}_{2,n}$ &$\hat{\gam}_{3,n}$& $\hat{\al}_{1,n}$&$\hat{\al}_{2,n}$\\\hline
        50 & 0.03 &20 & 1.938803 &0.971131 & -0.1962665
        &-0.3929058& 0.4998789 \\
        &&&(0.09929237)& (0.04901481) &(0.02746071) & (0.2725751) & (0.1483873) \\
        100&0.01&40 &1.937239 &0.9714114 &  -0.1945229& -0.4052143& 0.5019145\\
        &&& (0.06365245)& (0.03432344) & (0.01212126)& (0.1968092)&(0.09551109)\\
        200&0.004&75 &1.939557& 0.9712576&-0.1944697&-0.3951321&0.5010839\\
        &&&(0.04637029)&(0.02273138)&(0.006458429)&(0.1411913)&(0.07150592)\\\hline
    \end{tabular}
    \caption{The means of the estimators are reported, with standard deviations in parentheses.}
    \label{tab1}
\end{table}


\section{Proofs}\label{proofs-section}

\subsection{Preparation}

We first introduce the notations used in this section.
\begin{itemize}
\item For a matrix $M$, $[M]_{ij}$ denotes $(i,j)$ element of $M$.
\item $C$ denotes generic positive constant whose value may vary depending on context.
\item For positive sequences $\{a_n\}$ and $\{b_n\}$, $a_n\lesssim b_n$ means that 
$$\limsup_{n\to\infty} \frac{a_n}{b_n}<\infty.$$
\item For short, we write $\tilde{\Sig}_i^l(\gam)=V^l(\gam)|I_i^l|$, $\mcd=\text{diag}(\{|I_i|\}_{i=1}^{{\bf M}_n})$ and $\tilde{\Sig}_{i,j}^{1,2}(\gam)=(V^{1,2}(\gam)|I_i^1\cap I_j^2|)_{ij}$.
\item For a random variable ${\bf X}$, we denote 
$$E_\Pi[{\bf X}]=E[ {\bf X} \mid \{S^{n,l}_i\}_{0\leq i\leq {\bf M}_{l,n}, l \in \{1,2\}}].$$
\item Let $N(\cdot, \cdot)$ be the Poisson random measure associated with $J$, and let $\tilde{N}(\cdot, \cdot)$ be the compensated Poisson random measure.
\item Let $\mcv_t(\al)=\int_0^t A(\al)ds=A(\al)t$, and $X_t^d=\int_0^t C(\gamma_0)dJ_t=C(\gamma_0)J_t$.
\item We define $\chi_i^l(t)=1_{I_i^l}(t)$ and correspondingly, we write 
\begin{align*}
    \chi(t)=\big((\chi_i^1(t))_{1\leq i\leq {\bf M}_{1,n}}^\top, (\chi_i^2(t))_{1\leq i\leq {\bf M}_{2,n}}^\top\big)^\top.
\end{align*}
\end{itemize}
Then, a simple calculation shows that
\EQ{\label{chi-chi-t-eq} \bigg[\int_0^{T_n}\chi_{Cz}(t)\chi_{Cz}(t)^\top dt\bigg]_{i,j}
=[Cz]_{\psi(i)}[Cz]_{\psi(j)}|I_i\cap I_j|.}

Throughout this section, we often use the following decomposition of $S_n(\gam)$:
$$S_n(\gam)=\tmcd^{1/2}(\gam)\begin{pmatrix}I &\tilde{G}(\gam)\\ \tilde{G}^\top(\gam)& I \end{pmatrix}\tmcd^{1/2}(\gam),$$
where $$\tmcd(\gam)=\diag(\tilde{\Sig}_i^l(\gam)_{1\leq i\leq \bf{M}_n}), \ \ \tilde{G}(\gam)=\Bigg\{\frac{\tilde{\Sig}^{1,2}_{i,j}(\gam)}{\sqrt{\tilde{\Sig}^1_i}\sqrt{\tilde{\Sig}^2_{j}}}\Bigg\}_{1\leq i\leq {\bf M}_{1,n}, 1\leq j\leq {\bf M}_{2,n}}.$$ 

Let 
$$\bar{\rho}=\sup_\gamma \frac{|V^{1,2}(\gamma)|}{\sqrt{V^1(\gam)V^2(\gam)}}.$$

We prepare the following results in Ogihara and Yoshida~\cite{ogi-yos14} and Ogihara~\cite{ogi23}.

\begin{lemma}[Lemma 2 in Ogihara and Yoshida~\cite{ogi-yos14}, and Lemma 3.5 in~\cite{ogi23}]\label{ogi23-lem}
Assume (A2). Then,
$\bar{\rho}<1$ and there exists a positive constant $C$ such that
$$\| \mcd^{1/2}\PB^k S_n^{-1}(\gamma)\mcd^{1/2}\| \leq C (1-\bar{\rho})^{-k-1}, 
\quad \| \mcd^{-1/2}\PB^k S_n(\gamma)\mcd^{-1/2}\| \leq C$$
for $k\in \{0,1,2,3,4\}$ and any $\gamma\in \Theta_\gamma$. Moreover, 
\EQ{\label{G-ineq} \lVert G\rVert \vee \lVert G^\top \rVert \leq 1.}
\end{lemma}

\begin{lemma}\label{asym-norm-aux-lem}
Assume (A1) and (A2). Then, there exists a positive constant $C$ such that
$$
E_\Pi\bigg[\int_0^{T_n} \int\Big|\left(\Delta_{t-} X^{d}\right)^{\top} S_{n}^{-1}(\gam) \chi_{C z}(t)\Big|^{2} \nu(d z) d t \bigg]\leq C\M_n \AS
$$
for any $\gam \in \Theta_\gam$.
\end{lemma}

\begin{proof}
The proof of Lemma 3.5 in Ogihara~\cite{ogi23} yields
\EQ{\label{invSn-eq}
S_{n}^{-1}=\tmcd^{-1 / 2} \sum_{p=0}^{\infty}(-1)^{p}\left(\begin{array}{cc}
0 & \tilde{G} \\
\tilde{G}^{\top} & 0\\
\end{array}\right)^{p} \tmcd^{-1 / 2},
}
and hence, we have
$$
\left|\left[ S_{n}^{-1}\right]_{i j}\right| \leq C\left[\mathfrak{S}_{n}\right]_{i j}
$$
for $k=1,\ldots, p_\gamma$, where 
$$
\mathfrak{S}_n=\mcd^{-1/2}\sum_{p=0}^\infty (1+p)\bar{\rho}^{p}
\left(\begin{array}{cc}
0 & G \\
G^{\top} & 0\\
\end{array}\right)^{p}\mcd^{-1/2}.
$$
\begin{note}
    {\colorr Because all the terms are nonnegative.}
\end{note}

Together with (\ref{chi-chi-t-eq}), we have
$$
\begin{aligned}
& \int_0^{T_n}\int \left(\left(\Delta_{t-} X^{d}\right)^{\top}  S_{n}^{-1} \chi_{C z}(t)\right)^{2} \nu(d z) d t \\
&\quad =\int_0^{T_n}\int \sum_{i_1,i_2,j_1,j_2} [S_n^{-1}]_{i_1,j_1} [S_n^{-1}]_{i_2,j_2}
\Delta_{i_1,t-}X^d \Delta_{i_2,t-}X^d[\chi_{Cz}(t)]_{j_1}[\chi_{Cz}(t)]_{j_2}\nu(dz)dt \\
&\quad \leq C \sum_{i_1,i_2,j_1,j_2} [\mathfrak{S}_n]_{i_1,j_1}[\mathfrak{S}_n]_{i_2,j_2}
\sup_t|\Delta_{i_1,t-}X^d|\sup_t|\Delta_{i_2,t-}X^d|
\cdot \int |Cz|^2\nu(dz) \cdot |I_{j_1}\cap I_{j_2}| \\
&\quad \leq C\bigg(\sum_i[\mcd^{-1/2}]_{ii}^2\sup_t |\Delta_{i,t-}X^d|^2 \bigg)
\lVert \mcd^{1/2}\mathfrak{S}_n \mcd^{1/2} \rVert^2 \lVert \mcd^{-1/2}\{|I_{j_1}\cap I_{j_2}|\}_{j_1,j_2=1}^{\M_n} \mcd^{-1/2}\rVert \\
&\quad \leq C\sum_i \sup_t \frac{|\Delta_{i,t-}X^d|^2}{|I_i|}.
\end{aligned}
$$
Here we used that
$$\lVert \mcd^{-1/2}\{|I_{j_1}\cap I_{j_2}|\}_{j_1,j_2=1}^{\M_n} \mcd^{-1/2}\rVert =\bigg\lVert \MAT{cc}{I & G \\ G^\top & I}\bigg\rVert \leq 1+\lVert G \rVert \vee \lVert G^\top \rVert \leq 2,$$
and 
$$\lVert \mcd^{1/2}\mathfrak{S}_n \mcd^{1/2}\rVert \leq \SUMP (p+1)\bar{\rho}^{p} (\lVert G\rVert^p \vee \lVert G^\top \rVert^p)<\infty$$
by (\ref{G-ineq}).

Therefore, Burkholder's inequality yields
\EQNN{& E_\Pi \bigg[\int_0^{T_n}\int \left(\left(\Delta_{t-} X^{d}\right)^{\top}  S_{n}^{-1} \chi_{C z}(t)\right)^{2} \nu(d z) d t\bigg] \\
&\quad \leq C\sum_iE_\Pi\bigg[\sup_t \frac{|\Delta_{i,t-}X^d|^2}{|I_i|}\bigg] \\
&\quad \leq C\sum_i|I_i|^{-1}E_\Pi\bigg[\sup_t \bigg|\int_{I_i\cap [0,t)}[C(\gamma_0)dJ_s]_{\psi(i)}\bigg|^2\bigg] \\
&\quad \leq C\sum_i|I_i|^{-1}|I_i|=C\M_n \AS
}

\end{proof}

\subsection{Asymptotic behavior of $\ges$}
In this subsection, our goal is to give
the $\sqrt{nh_n}$-consistency of $\ges$.
We first show the consistency of $\ges$.
\begin{proposition}\label{consgam}
Under (A1)-(A7), we have
$$\ges \TOP \gam_0.$$
\end{proposition}

\begin{proof}
For any $\al\in\Theta_\al$, we have 
\begin{equation}\label{DV-ineq}
\big|\mcd^{-1/2}\Delta \mcv(\al)\big|^2\leq \sup_\al |A(\al)|^2 \sum_{i,l} \big|I_i^l\big|\leq Cnh_n.
\end{equation}
Combined with Lemma~\ref{ogi23-lem}, we have 
\begin{align*}
     \bigg|\big(\D \mcv\big)^\top S_n^{-1}(\gam)\D \mcv\bigg|\leq ||\mcd^{1/2} S_n^{-1}(\gam)\mcd^{1/2}|| |\mcd^{-1/2}\D \mcv|^2\leq C(1-\bar{\rho})^{-1}nh_n.
\end{align*}
The stationarity and independence of the increments of L\'{e}vy processes lead to the equation $E_\Pi[\D X^d(\D X^d)^\top]=S_n(\gam_0)$,
and thus Lemma~\ref{ogi23-lem} and (\ref{DV-ineq}) yield
\EQNN{E_\Pi\big[\big|(\D \mcv)^\top S_n^{-1}(\gam)\D X^d\big|^2\big]&=(\D \mcv)^\top S_n^{-1}(\gam)E_\Pi\big[\D X^d(\D X^d)^\top\big]S_n^{-1}(\gam)\D \mcv \\
&\leq C(1-\bar{\rho})^{-2}|\mcd^{-1/2}\D \mcv |^2 \\
&=O_p(nh_n).
}
Hence we obtain
\EQN{\label{H1n-eq}
    \mbbh_n^1(\gam)&=-\frac{1}{2}\D X^\top S_n^{-1}(\gam)\D X-\frac{1}{2}\log\det S_n(\gam)\\
    &=-\frac{1}{2}(\D X^d)^\top S_n^{-1}(\gam)\D X^d-\frac{1}{2}\log\det S_n(\gam)+\sqrt{n}e_n(\gam),
}
where $(e_n(\gamma))_{n=1}^\infty$ denotes a general sequence of random variables satisfying $\sup_\gam |e_n(\gam)|=o_p(1)$.

Let 
$$\mcr_n(\gamma)=(\D X^d)^\top S_n^{-1}(\gam)\D X^d-\text{tr}\big(S_n^{-1}(\gam)S_n(\gam_0)\big),$$
then it follows from It\^{o}'s formula that
\begin{align*}
    &\mcr_n(\gam)\\
    &=\sum_{i,j}\big[S_n^{-1}(\gam)\big]_{ij}\big(\D_i X^d \D_j X^d-[S_n(\gam_0)]_{ij}\big)\\
    &=\sum_{i,j}\big[S_n^{-1}(\gam)\big]_{ij}\Bigg\{2\int_{I_i}\D_{j,t-}X^d dX_t^{d,\psi(i)}+\int_{I_i\cap I_j}\int [C(\gam_0) z]_{\psi(i)}[C(\gam_0) z]_{\psi(j)} \tilde{N}(ds,dz)\Bigg\}\\
    &=2\int_0^{nh_n} \int \big(\D_{t-} X^d\big)^\top S_n^{-1}(\gam)\chi_{Cz}(t)d\tilde{N}(dt,dz)+\int_0^{nh_n}\int \big(\chi_{Cz}(s)\big)^\top S_n^{-1}(\gam)\chi_{Cz}(s) \tilde{N}(ds,dz).
\end{align*}
Applying It\^o isometry and Lemma~\ref{asym-norm-aux-lem}, we have 
\begin{align*}
&E_\Pi\Bigg[\Bigg|\int_0^{nh_n} \big(\D_{t-} X^d\big)^\top S_n^{-1}(\gam)\chi_{Cz}(t)d\tilde{N}(dt,dz)\Bigg|^2\Bigg] \\
&\quad = E_\Pi \bigg[\int_0^{nh_n} \Big| \big(\D_{t-} X^d\big)^\top S_n^{-1} \chi_{Cz}(t)\Big| \nu(dz)dt\bigg] \leq C\M_n.
\end{align*}
Moreover, we obtain
\begin{align*}
&E_\Pi\Bigg[\Bigg|\int_0^{nh_n}\int \big(\chi_{Cz}(s)\big)^\top S_n^{-1}(\gam)\chi_{Cz}(s) \tilde{N}(ds,dz)\Bigg|^2\Bigg]\\
    &\leq C\int_0^{nh_n}\Bigg|(\chi(s))^\top S_n^{-1}(\gam)\chi(s)\Bigg|^2ds\\
    &\leq C \int_0^{nh_n}\|\mcd^{1/2}S_n^{-1}(\gam)\mcd^{1/2}\|^2|\mcd^{-1/2}\chi(s)|^4ds\\
    &\leq C\int_0^{nh_n} \bigg(\sum_i\frac{[\chi(s)]_i}{|I_i|}\bigg)^2 ds.
\end{align*}
We obtain 
$$\int_0^{nh_n} \bigg(\sum_i\frac{[\chi(s)]_i}{|I_i|}\bigg)^2 ds
=\sum_{i,j}\frac{|I_i\cap I_j|}{|I_i||I_j|}\leq (1+\|G\|) \sum_i|I_i|^{-1}.$$
Together with Lemma~\ref{ogi23-lem}, we obtain
\EQ{\label{mcr-eq} E_\Pi[|\mcr_n(\gamma)|^2]\leq C\M_n + C\sum_{l=1,2} \sum_{i=1}^{\M_{l,n}} |I_i^l|^{-1} 
\leq C\M_n + C\M_n \underbar{r}_n^{-1}.}
Since (\ref{Mn-order-eq}) and (A3) yield
$$\frac{\M_n \underbar{r}_n^{-1}}{n^2}
=\frac{\M_n}{n}\cdot \frac{h_n\underbar{r}_n^{-1}}{nh_n} \TOP 0,$$
we have $\mcr_n(\gam)=o_p(n)$.
We also obtain 
\EQ{\label{mcr-eq3} \PB^k\mcr_n(\gam)=o_p(n)} for $k\in \{0,1,2,3\}$ by a similar argument.

Moreover, Lemma~\ref{ogi23-lem}
and the equation $E_\Pi[|\mcd^{-1/2}\Delta X^d|^2]={\rm tr}(\mcd^{-1}S_n(\gam_0))\leq C\M_n$ yield
\EQ{\label{mcr-eq2} \sup_{\gam\in \overline{\Theta_\gam}} |\PB^{k+1}\mcr_n(\gam)|
\leq \sup_{\gam\in \overline{\Theta_\gam}} \|\mcd^{1/2} \PB^{k+1} S_n^{-1}(\gam) \mcd^{1/2}\| 
\left(|\mcd^{-1/2}\Delta X^d|^2 + {\rm tr}(\mcd^{-1}S_n(\gam_0))\right) = O_p(n).
}

By compactness and convexity of $\overline{\Theta_\gamma}$, for every $\delta>0$ there exists a finite $\delta$-net. On the event that the Lipschitz constant of $n^{-1}\partial_\gamma^k\mcr_n$ is bounded by $K$, its supremum in absolute value over $\overline{\Theta_\gamma}$ is bounded by its maximum in absolute value over the net plus $K\delta$. For each fixed net, this maximum converges to zero in probability. By tightness of the Lipschitz constants, we obtain 
$$\sup_{\gam\in\Theta_\gam} |\PB^k\mcr_n(\gam)|=o_p(n)$$
by first choosing $K$ sufficiently large, then choosing $\delta$ sufficiently small, and finally letting $n\to\infty$.

Together with (\ref{H1n-eq}), we get
\begin{equation}
    \nn \mbbh_n^1(\gam)=-\frac{1}{2}\text{tr}\big(S_n^{-1}(\gam)S_n(\gam_0)\big)-\frac{1}{2}\log\det S_n(\gam)+ne_n(\gam),
\end{equation}
and consequently,
\begin{equation}
     \nn \mbbh_n^1(\gam)-\mbbh_n^1(\gam_0)=-\frac{1}{2}\text{tr}\big(S_n^{-1}(\gam)(S_n(\gam_0)-S_n(\gam)\big)-\frac{1}{2}\log\frac{\det S_n(\gam)}{\det S_n(\gam_0)}+ne_n(\gam).
\end{equation}

Similar estimates for the derivatives $\p_\gam^k(\mbbh_n^1(\gam)-\mbbh_n^1(\gam_0))$ give
\begin{equation}
    \label{sescale}\sup_\gam\Bigg|\frac{1}{n}\p_\gam^k(\mbbh_n^1(\gam)-\mbbh_n^1(\gam_0))+\frac{1}{2n}\p_\gam^k\text{tr}\big(S_n^{-1}(\gam)(S_n(\gam_0)-S_n(\gam)\big)+\frac{1}{2n}\p_\gam^k\log\frac{\det S_n(\gam)}{\det S_n(\gam_0)}\Bigg|=o_p(1),
\end{equation}
for $k=0,1,2,3$.
Then, we can obtain the consistency of $\HBN$ in the same way as the diffusion case \cite{ogi23} since $\text{tr}\big(S_n^{-1}(\gam)(S_n(\gam_0)-S_n(\gam)\big)$ and $\log\frac{\det S_n(\gam)}{\det S_n(\gam_0)}$ no longer depend on $J$.
\end{proof}

We note that \eqref{sescale} and the consistency of $\ges$ lead to 
\EQN{\label{convfig} \int_0^1\frac{1}{n}\PB^2 \mbbh^1_n(\gam_0+u(\ges-\gam_0))du + \mci_\gam =o_p(1),}
from the estimates in \cite{ogi23}.
We write $S_n=S_n(\gamma_0)$, $\tmcd=\tmcd(\gamma_0)$ and $\PB S_n=\PB S_n(\gamma_0)$ throughout this subsection.
Let $D_{l}=\operatorname{diag}\left(\left(\left|I_{i}^{l}\right|\right)_{i}\right)$ for $l\in\{1,2\}$.
Let
\EQNN{G_n^1 &=\bigg(\frac12(\Delta X^d)^\top (\partial_{\gamma_i} S_n^{-1}h_n)\Delta X^d +\frac{h_n}{2}\frac{\PB \det S_n}{\det S_n}\bigg)_{i=1}^{p_\gamma}.
}
Then, the estimates derived in the proof of Proposition \ref{consgam}, we have
\EQ{\label{G1-eq} -\frac{h_n}{\sqrt{T_n}}\PB \mbbh^1_n(\gam_0)
=\frac{1}{\sqrt{T_n}}G_n^1+o_p(1) 
=\frac{1}{2\sqrt{T_n}}((\Delta X^d)^\top (\PB S_n^{-1}h_n)\Delta X^d -h_n\text{tr}(\PB S_n^{-1}S_n))+o_p(1).}
To verify the behavior of the one-step estimator $\aes$, we now show the $\sqrt{nh_n}$-consistency of $\ges$.
For this purpose, we collect some auxiliary lemmas.

\begin{lemma}\label{asym-norm-aux-lem2}
Assume (A1), (A2) and (A4). Then, 
    $$T_n^{-1/2}G_n^1 = \frac{h_n}{2\sqrt{T_n}}\int_0^{T_n}\int (\chi_{Cz}(t))^\top \PB S_n^{-1} \chi_{Cz}(t)\tilde{N}(dt,dz) + o_p(1).$$
\end{lemma}

\begin{proof}

It\^o's formula for L\'evy processes (Theorem 4.5.7 in Chapter I of Jacod and Shiryaev~\cite{jac-shi03}) shows that
\begin{equation}
\label{XdInvSXd-eq2}
\begin{split}
 &\left(  \Delta X^d\right)^{\top}   \partial_\gamma S_n^{-1} \Delta X^d \\ 
 &  =  2 \int_0^{T_n}\int \left(\Delta_{t-} X^d\right)^{\top}   \partial_\gamma S_n^{-1} \chi_{{Cz}}(t) \tilde{N}(dt, dz) +  \int_0^{T_n}\int \left(  \chi_{{Cz}}(t)\right)^{\top}   \partial_\gamma S_n^{-1} \chi_{{Cz}}(t) N(dt, dz). 
 \end{split}
\end{equation}

Burkholder's inequality, (\ref{Mn-order-eq}) and a similar estimate to Lemma~\ref{asym-norm-aux-lem} yield
\begin{equation}
\begin{split}
 &E_{ \Pi } \left[\left|\int_0^{T_n}\int \left(\Delta_{t-} X^d\right)^{\top}   \partial_\gamma S_n^{-1} \chi_{{Cz}}(t) \tilde{N}(dt, dz)\right|^{ 2 }\right]\\ 
 &  \leq  C E_{ \Pi } \left[\left(\int_0^{T_n} \int_{  } (\left(\Delta_{t-} X^d\right)^{\top}   \partial_\gamma S_n^{-1} \chi_{{Cz}}(t))^{ 2 } \nu (dz) dt\right)\right]\\ 
 &  =   O_p(n). 
 \end{split}
\end{equation}
Then, Lemma B.1 in~\cite{ogi-ueh23} yields
\begin{equation}
\label{XdInvSXd-eq3}
\int_0^{T_n}\int \left(\Delta_{t-} X^d\right)^{\top}   \partial_\gamma S_n^{-1} \chi_{{Cz}}(t) \tilde{N}(dt, dz) =   O_p(n^{  1/2 }). 
 \end{equation}
Therefore, (\ref{XdInvSXd-eq2}) and (\ref{XdInvSXd-eq3}) yield
\begin{equation}
\label{XdInvSXd-eq}
\begin{split}
 &\left(  \Delta X^d\right)^{\top}   \partial_\gamma S_n^{-1} \Delta X^d\\ 
 &  =  \int_0^{T_n}\int \left(  \chi_{{Cz}}(t)\right)^{\top}   \partial_\gamma S_n^{-1}\chi_{{Cz}}(t) \tilde{N}(dt, dz) +  \int_0^{T_n} \int_{  } \left(  \chi_{{Cz}}(t)\right)^{\top}   \partial_\gamma S_n^{-1}   \chi_{{Cz}}(t) \nu (dz) dt +    O_p( n^{  1/2 }). 
 \end{split}
\end{equation}
Since (\ref{chi-chi-t-eq}) yields
\begin{equation}
\begin{split}
 &\int_0^{T_n} \int_{  } \left(  \chi_{{Cz}}(t)\right)^{\top}   \partial_\gamma S_n^{-1}   \chi_{{Cz}}(t) \nu (dz) dt\\ 
 &  = \int_0^{T_n} \int_{  }  \text{tr}(     \partial_\gamma S_n^{-1}     \chi_{{Cz}}(t)(  \chi_{{Cz}}(t))^\top) \nu (dz) dt\\ 
 &  =  \text{tr}\bigg(     \partial_\gamma S_n^{-1}   \int_0^{T_n}   \int_{  }     \chi_{{Cz}}(t)(  \chi_{{Cz}}(t))^\top \nu(dz) dt\bigg)\\ 
 &  =  \text{tr}(     \partial_\gamma S_n^{-1}   S_n), 
 \end{split}
\end{equation}
together with (\ref{XdInvSXd-eq}), we have
\begin{equation}
\begin{split}
 &  h_n \left(  \Delta X^d\right)^{\top}   \partial_\gamma S_n^{-1} \Delta X^d \\ 
 &  =   h_n \int_0^{T_n}\int \left(  \chi_{{Cz}}(t)\right)^{\top}   \partial_\gamma S_n^{-1}   \chi_{{Cz}}(t) \tilde{N}(dt, dz) +   h_n  \text{tr}(     \partial_\gamma S_n^{-1}   S_n) +    o_p( n^{  1/2 }  h_n^{  1/2 }). 
 \end{split}
\end{equation}
\end{proof}

Since Assumption \ref{convav} ensures $T_n^{-1/2}G^1_n=O_p(1)$, we obtain  
$$\sqrt{nh_n}(\ges-\gam_0)=O_p(1),$$
by Taylor expansion and \eqref{convfig}. 


\subsection{Asymptotic behavior of $\tilde{\al}_n$ and $\aes$}

First we provide the $\sqrt{nh_n}$-consistency of the initial estimator $\tilde{\al}_n$.

\begin{lemma}
Assume (A1), (A2), (A3) and (A7). Then, we have $\sqrt{nh_n}(\tilde{\al}_n-\al_0)=O_p(1).$
\end{lemma}

\begin{proof}
Let $\mbbg_2(\al)=||A(\al_0)-A(\al)||^2$.
We have
\begin{align*}
    &\frac{1}{nh_n}(\mbbg_{2,n}(\al)-\mbbg_{2,n}(\al_0))\\
    &=\frac{1}{nh_n}\sum_{i=1}^{{\bf M}_{1,n}}[A(\al_0)-A(\al)]_1^2|I_i^1| + \frac{1}{nh_n}\sum_{j=1}^{{\bf M}_{2,n}} [A(\al_0)-A(\al)]_2^2|I_j^2|\\
    &+\frac{2}{nh_n}\sum_{i=1}^{{\bf M}_{1,n}} [A(\al_0)-A(\al)]_1\D^1_i X^d+\frac{2}{nh_n}\sum_{j=1}^{{\bf M}_{2,n}} [A(\al_0)-A(\al)]_2\D^2_j X^d\\
    &=\mbbg_2(\al)+o_p((nh_n)^{-1}) +\frac{2}{nh_n}\sum_{i=1}^{{\bf M}_{1,n}} [A(\al_0)-A(\al)]_1\D^1_i X^d+\frac{2}{nh_n}\sum_{j=1}^{{\bf M}_{2,n}} [A(\al_0)-A(\al)]_2\D^2_j X^d,
\end{align*}
and it follows from the mean value theorem and Jensen's inequality that for all $\del>0$, 
\begin{align*}
    &\sqrt{nh_n}E\Bigg[\sup_{|\al-\al_0|<\del,\al\in\Theta_\al}\Bigg|\frac{1}{nh_n}(\mbbg_{2,n}(\al)-\mbbg_{2,n}(\al_0))-\mbbg_2(\al)\Bigg|\Bigg]\\
    &\lesssim \frac{1}{\sqrt{nh_n}}\sup_{|\al-\al_0|<\del,\al\in\Theta_\al}|A(\al_0)-A(\al)| \sqrt{E[|X^d_{nh_n}|^2]}\\
    &\lesssim \delta.
\end{align*}
An easy application of \cite[Theorem 3.2.5]{VanWel96} leads to 
$$\sqrt{nh_n}(\tilde{\al}_n-\al_0)=O_p(1).$$
\end{proof}

\subsection{Stochastic expansion of one-step estimator}
By the definition of $\aes$, the Taylor expansion implies that
\begin{align*}
    &\sqrt{nh_n}(\aes-\al_0)\\
    &=-\Bigg(\frac{1}{nh_n}\p_\al^2\mbbh_n^2(\tilde{\al}_n)\Bigg)^{-1}\frac{1}{\sqrt{nh_n}}\p_\al \mbbh_n^2(\al_0)\\
    &+\Bigg[I_{p_\al}-\Bigg(\frac{1}{nh_n}\p_\al^2\mbbh_n^2(\tilde{\al}_n)\Bigg)^{-1}\frac{1}{nh_n}\int_0^1\p_\al^2\mbbh_n^2(\al_0+s(\tilde{\alpha}_n-\al_0))ds \Bigg] \sqrt{nh_n}(\tilde{\al}_n-\al_0),
\end{align*}
Once, assume that for small enough $\del>0$,
\EQN{\label{td}\frac{1}{(nh_n)^{3/2}}\sup_{|\al-\al_0|<\del, \al\in\Theta_\al}|\p_\al^3\mbbh_n^2(\al)|=o_p(1),}
and that
\EQN{\label{conv:sd}
\frac{1}{nh_n}\p_\al^2\mbbh_n^2(\al_0)+\mathcal{I}_\al=o_p(1).
}
Then it follows from $\sqrt{nh_n}$-consistency of $\tilde{\al}_n$ that for all $\ep,\ep'>0$, there exist $N\in\mbbn$ and $M$ such that for any $n\geq N$,
\begin{align*}
    &P\Bigg(\Bigg|\frac{1}{nh_n}\p_\al^2\mbbh_n^2(\tilde{\al}_n)+\mathcal{I}_\al\Bigg|>\ep\Bigg)\\
    &\leq P\Bigg(\Bigg|\frac{1}{(nh_n)^{3/2}}\int_0^1\p_\al^3\mbbh_n^2(\al_0+s(\tilde{\alpha}_n-\al_0))ds\Bigg|\cdot|\sqrt{nh_n}(\tilde{\al}_n-\al_0)|>\frac{\ep}{2}\Bigg)+P\Bigg(\Bigg|\frac{1}{nh_n}\p_\al^2\mbbh_n^2(\al_0)+\mathcal{I}_\al\Bigg|>\frac{\ep}{2}\Bigg)\\
    &\leq P\Bigg(\frac{1}{(nh_n)^{3/2}}\sup_{|\al-\al_0|<\del, \al\in\Theta_\al}|\p_\al^3\mbbh_n^2(\al)|>\frac{\ep}{2M}\Bigg)+P(|\sqrt{nh_n}(\tilde{\al}_n-\al_0)|>M) \\
    &\quad +P\Bigg(\Bigg|\frac{1}{nh_n}\p_\al^2\mbbh_n^2(\al_0)+\mathcal{I}_\al\Bigg|>\frac{\ep}{2}\Bigg)\\
    &\leq \ep'.
\end{align*}
Similarly, we also have 
$$\frac{1}{nh_n}\int_0^1\p_\al^2\mbbh_n^2(\al_0+s(\tilde{\alpha}_n-\al_0))ds+\mathcal{I}_\al=o_p(1),$$
and Slutsky's lemma gives
\EQ{\label{alpha-est-eq} \sqrt{nh_n}(\aes-\al_0)=-\Bigg(\frac{1}{nh_n}\p_\al^2\mbbh_n^2(\tilde{\al}_n)\Bigg)^{-1}\frac{1}{\sqrt{nh_n}}\p_\al \mbbh_n^2(\al_0)+o_p(1).}
By mimicking the estimates in \cite[(3.47)]{ogi23}, for all natural number  $l\leq \tilde{m}$, we have
\begin{align*}
    \Bigg|\p_\al \mcv(\alpha_0)^\top \int_0^1 \p_\gam^{l}S_n^{-1}(\gam_0+s(\ges-\gam_0))ds[(\ges-\gam_0)^{\otimes l}]\Delta X^d\Bigg|=O_p\bigg(n^{(2-l)/2}h_n^{(1-l)/2}\bigg),
\end{align*}
and the boundedness of the derivative of $A(\cdot)$ and the estimates in \cite[(3.49)]{ogi23} imply that
\begin{align*}
     |\p_\al \mcv(\alpha_0)^\top \p_\gam^l S_n^{-1}(\gam_0)[(\ges-\gam_0)^{\otimes l}]\Delta X^d|\leq |\ges-\gam_0|^l|\p_\al \mcv(\alpha_0)^\top \p_\gam^l S_n^{-1}(\gam_0)\Delta X^d|=O_p((nh_n)^{\frac{1-l}{2}}).
\end{align*}
Hence, Taylor's expansion leads to
\begin{equation}
    \label{alscore} \frac{1}{\sqrt{nh_n}}\p_\al \mbbh_n^2(\al_0)=\frac{1}{\sqrt{nh_n}}G_n^2+o_p(1),
\end{equation}
where 
$$G_n^2= \bigg(\partial_{\alpha_i} \Delta \mcv(\alpha_0)^\top S_n^{-1}\Delta X^d\bigg)_{i=1}^{p_\alpha}.$$
From similar estimates and Taylor's expansion, we obtain \eqref{td} and
\begin{equation}\label{sedrift}
    \frac{1}{nh_n}\p_\al^2\mbbh_n^2(\al_0)=-\frac{1}{nh_n}\p_\al \mcv(\alpha_0)^\top S_n^{-1}\p_\al \mcv(\alpha_0)+o_p(1).
\end{equation}
Then, since the first term of the right-hand-side no longer depends on the driving noise $J$, \eqref{conv:sd} follows from \cite{ogi23}.

\bigskip

\subsection{Asymptotic normality}\label{asym-norm-subsection}

Let $\mcg_k=\mcf_{\check{s}_k}\vee \sigma(\Pi_n)$.
Let 
$$\mcl_n=[T_nn^{-\beta}], \quad \beta = \min\{\delta /2, 1/(2m)\},$$
where $m$ defined in (\ref{bcond}) and $\delta$ defined in (A3). Then, we obtain  
\EQ{\label{mcl-eq} \mcl_{n} \rightarrow \infty, \quad 
\mcl_n^{-1}(h_n\underbar{r}_n^{-1})^4\TOP 0,
\quad 
\mcl_{n} n^{\beta /2}T_n^{-1}\to 0}
as $n \rightarrow \infty$.
Let $\check{s}_k=\frac{k T_{n}}{\mcl_{n}}$ for $0 \leq k \leq \mcl_{n}$. Let $K^{k}=\left(\check{s}_{k-1}, \check{s}_{k}\right]$,
and let $\xi_k=(\xi_{k,l})_{l=1}^{p_\alpha}$ and $\zeta_k=(\zeta_{k,l})_{l=1}^{p_\gamma}$, where
for each $l$, $\alpha_l$ and $\gamma_l$ denote the $l$-th elements of $\alpha$ and $\gamma$, respectively, 
$$
\zeta_{k,l}=-\frac{h_n}{2\sqrt{T_{n}}} \int_{K^k} \int \chi_{Cz}(t)^{\top} \partial_{\gamma_l} S_{n}^{-1} \chi_{Cz}(t)  \tilde{N}(dt, dz),
$$
and
$$
\xi_{k,l}=\frac{1}{\sqrt{T_{n}}} \int_{K^k} \int \partial_{\alpha_l}\Delta \mcv(\alpha_0)^{\top} S_{n}^{-1} \chi_{Cz}(t)  \tilde{N}(dt, dz).
$$

Then, we obtain
\EQ{\label{Gn-approx-eq} -\frac{1}{\sqrt{T_{n}}} G_n^1=\sum_{k} \zeta_{k}+o_{p}(1),
\quad \frac{1}{\sqrt{T_{n}}} G_n^2=\sum_{k} \xi_{k}}
by Lemma~\ref{asym-norm-aux-lem2} and the fact that
\EQNN{\frac{1}{\sqrt{T_n}}G_n^2 
=\frac{1}{\sqrt{T_n}}\PA \Delta\mcv(\alpha_0)^\top S_n^{-1}\Delta X^d 
=\frac{1}{\sqrt{T_n}}\int_0^{T_n}\int 
\PA \Delta\mcv(\alpha_0)^\top S_n^{-1}\chi_{Cz}(t)\tilde{N}(dt,dz).
}

\begin{note}
{\colorr 

$$E[\xi_n^2|\mcf_n]
=\frac{1}{T_n} \int_0^{T_n}\int (\PA \Delta \mcv(\alpha)S_n^{-1}C\chi_{Cz}(t))^2dt\nu (dz)+O_p(T_n^{-1}T_n^2h_n^2).
$$
$$E[\xi_n\zeta_n|\mcf_n]
=\frac{1}{T_n} \int_0^{T_n}\int \big\{\PA \Delta \mcv(\alpha)S_n^{-1}\chi_{Cz}(t)\chi_{Cz}(t)^\top \PB S_n^{-1}\chi_{Cz}(t)\big\}dt\nu (dz).
$$

Let $\theta_{u}=u \hat{\theta}_n+(1-u) \theta_{0}$.

$$
\begin{aligned}
 -\partial_{\theta}Q_{n}\left(\theta_{0}\right)
&=\partial_{\theta} Q_{n}(\hat{\theta}_{n})
-\partial_{\theta} Q_n(\theta_{0}) \\
& =\int_{0}^{1} \partial_{\theta}^{2} Q_n\left(\theta_{u}\right) d u(\hat{\theta}_n-\theta_{0})  \\
& =\left(\begin{array}{c}
\int_{0}^{1} \partial_{\gamma}^{2} Q_n(\theta_u) du \cdot (\hat{\gamma}_{n}-\gamma_{0})+\int_{0}^{1} \partial_{\gamma}\partial_{\alpha} Q_n(\theta_u) du \cdot (\hat{\alpha}_{n}-\alpha_{0}) \\
\int_{0}^{1} \partial_{\alpha}^{2} Q_n(\theta_u) du \cdot (\hat{\alpha}_{n}-\alpha_{0})+\int_{0}^{1} \partial_{\alpha}\partial_{\gamma} Q_n(\theta_u) du \cdot (\hat{\gamma}_{n}-\gamma_{0})
\end{array}\right).
\end{aligned}
$$

$$
\frac{1}{\sqrt{T_n}} G_{n}^{\gamma} \approx -\frac{h_{n}}{\sqrt{T_{n}}} \PB Q_n, \quad \frac{1}{\sqrt{T_{n}}} G_{n}^{\alpha}=\frac{1}{2 \sqrt{T_{n}}} \partial_{\alpha} Q_{n}.
$$

$$
\begin{aligned}
E[\xi_{n}^{2} \mid \mcf_{n}]
& =\frac{1}{T_{u}} \int_{0}^{T_{u}} \int\left(\partial_{\alpha} \Delta \mcv(\alpha) S_n^{-1} \chi_{Cz}(t)\right)^{2} \nu(d z) \\
& =\frac{1}{T_n} \int \partial_{\alpha} \Delta \mcv(\alpha)^{\top} S_n^{-1} \int_{0}^{T_n} \chi_{C z}(t) \chi_{C z}(t)^{\top} d t S_n^{-1} \partial_{\alpha} \Delta \mcv(\alpha) \nu(d z) \\
& =\frac{1}{T_{n}} \partial_{\alpha} \Delta \mcv^{\top}  S_{n}^{-1} \partial_\alpha \Delta \mcv.
\end{aligned}
$$

$$
\begin{aligned}
 E\left[\xi_{n} \zeta_{n} \mid \mcf_{n}\right]
&= \frac{h_{n}}{T_{n}} \int_{0}^{T_{n}} \int \partial_\alpha \Delta \mcv S_{n}^{-1} \chi_{C z}(t) \chi_{C z}(t)^{\top} \partial_\alpha S_{n}^{-1} \chi_{C z}(t) d t \nu(d z) \\
&= \frac{h_{n}}{T_{n}} \sum_{i}\left|I_{i}\right|\left[S_{u}^{-1} \partial \alpha \Delta \mcv\right]_{i} \int[C z]_{\psi_{(i)}} \mcu_{i,i} \nu(d z) \\
&\quad +\frac{h_{n}}{T_{n}} \sum_{i \neq i}\left|I_{i} \cap I_{i'}\right|\left[S_{n}^{-1} \partial_{\alpha} \Delta \mcv\right]_{i} \int [Cz]_{\psi(i)}\left(\mcu_{i',i'}+2 \mcu_{i,i'}\right) \nu(d z).
\end{aligned}
$$
}
\end{note}

\begin{lemma}\label{invSdelV-est-lemma}
Assume (A2). There exists a positive constant $C$ such that
$$
\left|\left[S_{n}^{-1} \PA \Delta \mcv(\alpha)\right]_i\right| \leq C
$$
for any $\alpha$ and $1 \leq i \leq M_{n}$.
\end{lemma}

\begin{proof}
Let 
\EQ{\label{S0-def} S_{0}=\MAT{ll}{I & \tilde{G} \\ \tilde{G}^\top & I}.}
Then, since 
$$S_n^{-1}=(\tmcd^{1/2}S_0\tmcd^{1/2})^{-1}
=\left\{\tmcd(\tmcd^{-1/2}S_0\tmcd^{-1/2})\tmcd\right\}^{-1}
=\tmcd^{-1}(\tmcd^{1/2}S_0^{-1}\tmcd^{1/2})\tmcd^{-1},$$ 
we obtain
\EQN{\label{invSdelV-est-lemma-eq1}
\bigg|\left[S_{n}^{-1} \PA \Delta \mcv(\alpha)\right]_i\bigg|
&=\bigg|[\tmcd^{-1}]_{ii} \sum_{j}\left[\tmcd^{1/2} S_{0}^{-1} \tmcd^{1/2}\right]_{ij} [\tmcd^{-1}]_{jj}\left[\PA \Delta \mcv\right]_{j} \bigg|, \\
& \leq C\sum_{j} \bigg|\frac{[\mcd^{1/2} S_{0}^{-1} \mcd^{1/2}]_{i j} \PA A^{\psi(j)}(\alpha)\left|I_{j}\right|}{V^{\psi(i)}\left|I_{i}\right| V^{\psi(j)}\left|I_{j}\right|}\bigg| \\
& \leq C \sum_{j} \frac{|[\mcd^{1/2} S_0^{-1} \mcd^{1/2}]_{i j}|}{\left|I_{i}\right|}.
}
Since all elements of $G$ is nonnegative,
the equation $(I+A)^{-1}=\SUMP (-1)^pA^p$ for $\lVert A\rVert <1$  yields
\EQNN{&\sum_j|[\mcd^{1/2}S_0^{-1} \mcd^{1/2}]_{ij}| \\
&\quad = \sum_j \bigg|\bigg[\mcd^{1/2}\SUMP (-1)^p
\MAT{cc}{0 & \tilde{G} \\ \tilde{G}^\top & 0}^p \mcd^{1/2}\bigg]_{ij}\bigg| \\
& \quad \leq \SUMP \bar{\rho}^p \sum_j \bigg|\bigg[\mcd^{1/2}\MAT{cc}{0 & G \\ G^\top & 0}^p \mcd^{1/2}\bigg]_{ij}\bigg|.
}

When $p=2k+1$ for some positive integer $k$, then 
$$\bigg[\mcd^{1/2} (-1)^p
\MAT{cc}{0 & \tilde{G} \\ \tilde{G}^\top & 0}^p \mcd^{1/2}\bigg]_{ij}=0$$
if $i\vee j\leq {\bf M}_{1,n}$ or $i\wedge j>{\bf M}_{1,n}$.
If $i>{\bf M}_{1,n}$ and $j\leq {\bf M}_{1,n}$, then by setting $i'=i-{\bf M}_{1,n}$ and $j_0=i'$, we have
\EQQ{\bigg[\mcd^{1/2}\MAT{cc}{0 & G \\ G^\top & 0}^p \mcd^{1/2}\bigg]_{ij} =\sum_{\substack{i_1,\cdots, i_k \\ j_1,\cdots, j_k}}\sqrt{|I_{i'}^2|}\bigg\{\prod_{l=1}^k
\bigg(\frac{|I_{j_{l-1}}^2\cap I_{i_l}^1|}{\sqrt{|I_{j_{l-1}}^2|}\sqrt{|I_{i_l}^1|}}
\frac{|I_{i_l}^1\cap I_{j_l}^2|}{\sqrt{|I_{i_l}^1|} \sqrt{|I_{j_l}^2|}} \bigg) \bigg\}\frac{|I_{j_k}^2\cap I_{j}^1|}{\sqrt{|I_{j_k}^2|}}.
}
Then, repeated use of the equations
$$
\begin{aligned}
& \sum_{j} \frac{\left|I_{j_k}^{2} \cap I_{j}^{1}\right|}{\sqrt{\left|I_{j_k}^{2}\right|}}\leq \sqrt{\left|I_{j _k}^{2}\right|}, 
\quad \sum_{j_l} \frac{\left|I_{i_l}^{1} \cap I_{j_l}^{2}\right|}{\sqrt{\mid I_{i_l}^{1}}\left|\sqrt{\mid I_{j_l}^{2}}\right|} \sqrt{\left|I_{j_l}^{2}\right|}\leq \sqrt{\left|I_{i_l}^{1}\right|}, \\
& \sum_{i_l} \frac{\left|I_{j_{l-1}}^{2} \cap I_{i_l}^{1}\right|}{\sqrt{|I_{j_{l-1}}^{2}|} \sqrt{\left|I_{i_l}^{1}\right|}} \sqrt{\left|I_{i_l}^{1}\right|}\leq \sqrt{\left|I_{j_{l-1}}^{2}\right|}
\end{aligned}
$$
yields
$$
\sum_j\left[\mcd^{1/2}\MAT{cc}{0 & G \\ G^\top & 0}^{p} \mcd^{1/2}\right]_{i j}\leq \left|I_{i}^{2}\right| .
$$
Together with a similar argument for other cases of $p$, we obtain
$$
\sum_j\left[\mcd^{1/2}\MAT{cc}{0 & G \\ G^\top & 0}^{p} \mcd^{1/2}\right]_{i j} \leq \left|I_{i}\right|
$$
for any $p \in \mathbb{Z}+$.

Then we have $\sum_j[\mcd^{1/2}S_0^{-1} \mcd^{1/2}]_{ij}|\leq C|I_i|$, and together with (\ref{invSdelV-est-lemma-eq1}), we obtain the conclusion.
\end{proof}

\begin{lemma}\label{zeta4-conv}
Assume (A1), (A2), (A3), and (A6). Then,
\EQQ{\sum_kE[|\zeta_k|^4\mid \FCSK]\TOP 0, \quad \sum_kE[|\xi_k|^4\mid \FCSK]\TOP 0}
as $n\to\infty$.
\end{lemma}

\begin{proof}
Kunita's first inequality (Theorem 4.4.23 in Applebaum) yields
$$
\begin{aligned}
& \sum_k E\left[|\xi_{k}|^{4} \mid \FCSK \right] \\
& \quad =\frac{C}{T_{n}^{2}} \sum_{k} E\bigg[\bigg|\int_{K^k} \int \PA \Delta \mcv^{\top} S_{n}^{-1} \chi_{C z}(t) \tilde{N}(dt dz)\bigg|^{4} ~\bigg|~ \FCSK \bigg] \\
& \quad \leq \frac{C}{T_{n}^{2}} \sum_k \bigg\{\bigg(\int_{K^k} \int\Big|\PA \Delta \mcv^{\top} S_{n}^{-1} \chi_{C z}(t)\Big|^{2} \nu(d z) dt\bigg)^{2} \\
& \qquad \qquad +\int_{K^{k}} \int \Big|\PA \Delta \mcv^{\top} S_{n}^{-1} \chi_{C z}(t)\Big|^{4} \nu(d z) d t\bigg\}.
\end{aligned}
$$

By using Lemma~\ref{invSdelV-est-lemma}, we have
$$
\begin{aligned}
&\sum_{k} E\left[|\xi_{k}|^{4} \mid \FCSK \right] \\
& \quad \leq \frac{C}{T_{n}^{2}} \sum_{k}
\bigg\{\bigg(\int_{K^k} \int \sum_{i, j}|[\chi_{Cz}(t)]_{i} [\chi_{Cz}(t)]_{j}| \nu(d z) d t\bigg)^{2} +\int_{K^k} \int \sum_{i_1,\cdots, i_4} \prod_{l=1}^{4}\left|\left[\chi_{Cz}(t)\right]_{i_l}\right| \nu(d z) d t\bigg\} \\
& \quad \leq \frac{C}{T_{n}^{2}} \sum_{k}\bigg\{\bigg(\sum_{i, j}\left|I_{i} \cap I_{j} \cap K^k\right|\bigg)^{2}+\sum_{i, j}\left|I_{i} \cap I_{j} \cap K^k\right|\bigg\} \\
& \quad \leq \frac{C}{T_{n}^{2}} \sum_{k}\left|K^k\right|^{2}+\frac{1}{T_{n}}=\frac{1}{T_{n}^{2}} \cdot \frac{T_{n}}{\mcl_{n}} \cdot T_{n}+\frac{1}{T_{n}} \rightarrow 0
\end{aligned}
$$
as $n \rightarrow \infty$.

Again, Kunita's first inequality yields
$$
\begin{aligned}
& \sum_{k} E\left[|\zeta_k|^{4} \mid \FCSK \right] \\
&\quad =\frac{h_{n}^{4}}{16T_{n}^2} \sum_k E\bigg[ \bigg| \int_{K^k} \int \chi_{Cz}(t)^{\top} \PB S_{n}^{-1} \chi_{Cz}(t) \tilde{N}(dt d z)\bigg|^{4}\bigg| \mcg_{k-1} \bigg] \\
& \quad \leq \frac{C h_{n}^{4}}{T_{n}^2} \sum_{k}
\bigg\{\bigg(\int_{K^k} \int \Big|\chi_{Cz}(t)^{\top} \PB S_n^{-1} \chi_{Cz}(t)\Big|^{2} \nu\left(dz\right) d t\bigg)^{2}
+\int_{K^k} \int\Big|\chi_{Cz}(t)^{\top} \PB S_n^{-1} \chi_{Cz}(t)\Big|^{4} \nu\left(dz\right) d t\bigg\} \\
&\quad =: \Lambda_1 + \Lambda_2.
\end{aligned}
$$
By Lemma~\ref{ogi23-lem}, for any $k$, we obtain
\EQ{\label{invSn-est} |[\partial_{\gam_k} S_n^{-1}]_{ij}|
\leq \frac{|[D^{1/2} \partial_{\gam_k} S_n^{-1} D^{1/2}]_{ij}|}{|I_i|^{1/2}|I_j|^{1/2}}
\leq \frac{\| D^{1/2} \partial_{\gam_k} S_n^{-1} D^{1/2}\|}{|I_i|^{1/2}|I_j|^{1/2}}
\leq C\underbar{r}_n^{-1}.}
$\int_{K^{k}} \prod_{j=1}^{4} [\chi_{Cz}(t)]_{i_j} d t$ is non-zero only if $i_1=\cdots = i_4$ or the case that $\{i_1,\cdots,i_4\}$ is equal to a set of indicators $\{i,i'\}$ with $i\leq {\bf M}_{1,n}<i'$ and $I_i\cap I_{i'}\neq \emptyset$.
Therefore, by setting $\mfa_{i,i'}(z)=|[Cz]_{\psi(i)}[Cz]_{\psi(i')}[\PB S_n^{-1}]_{i,i'}|$, we obtain
\EQNN{&\int_{K^k} \int \Big|\chi_{Cz}(t)^{\top} \PB S_n^{-1} \chi_{Cz}(t)\Big|^{2} \nu\left(dz\right) d t \\
& \quad \leq \sum_{i_1, \cdots i_4}\Big|\left[\PB S_n^{-1}\right]_{i_1, i_2}\left[\PB S_n^{-1}\right]_{i_3,i_4}\Big|\int_{K^{k}}\int  \prod_{j=1}^{4} |[\chi_{Cz}(t)]_{i_j}| \nu (dz)dt \\
&\quad \leq \sum_i|I_i\cap K^k|\int \mfa_{ii}^2\nu(dz)
+\SUMIIp |I_i\cap I_{i'}\cap K^k|\int \left\{4\mfa_{i,i'}\mfa_{i,i}+4\mfa_{i,i'}\mfa_{i',i'}+4\mfa_{i,i'}^2+2\mfa_{i,i}\mfa_{i',i'}\right\}\nu(dz) \\
&\quad \leq C\underbar{r}_n^{-2}\sum_i|I_i\cap K^k|
+C\underbar{r}_n^{-2}\SUMIIp |I_i\cap I_{i'}\cap K^k| \\
& \quad \leq C\underbar{r}_n^{-2}|K^k|.
}

Therefore, (\ref{mcl-eq}) yields
$$\Lambda_1 \leq \frac{Ch_n^4\underbar{r}_n^{-4}}{T_n^2}\sum_k|K^k|^2\leq \frac{Ch_n^4\underbar{r}_n^{-4}}{\mcl_n}\TOP 0$$
as $n\to\infty$.

Moreover, since $\int_{K^k}\prod_{k=1}^8|[\chi_{Cz}(t)]_{i_k}| dt$ is non-zero only if $i_1=\cdots = i_8$ or the case that $\{i_1,\cdots, i_8\}$ is equal to a set of indicators $\{i,j\}$ with $i\leq {\bf M}_{1,n}<j$ and $I_i\cap I_j\neq \emptyset$.
Then, together with (\ref{invSn-est}) and (A3), we have
$$\Lambda_2\leq \frac{Ch_n^4}{T_n^2}
\sum_k\bigg\{\underbar{r}_n^{-4}\sum_i|I_i\cap K^k| + \underbar{r}_n^{-4} \SUMIIp |I_i\cap I_{i'}\cap K^k|\bigg\}
\leq \frac{Ch_n^4\underbar{r}_n^{-4}}{T_n^2}\cdot T_n\TOP 0$$
as $n\to\infty$, which completes the proof.

\end{proof}

\noindent
{\bf Proof of Theorem~\ref{asym-norm-thm}.}

We have that the matrix $\Gamma_1$ is nonnegative definite because $F_n\TOP \Gamma_1$ by (A6) and
\[
u^\top F_n u
=
\frac{1}{T_n}
 \int_0^{T_n}\int
 \bigl(u^\top\mathcal{S}(t,z)\bigr)^2\nu(dz)dt
\geq 0
\]
for every $u\in\mathbb{R}^{p_\gamma+p_\alpha}$.

We have
\EQNN{h_n\PB \mbbh_n^1(\gamma_0) 
&= h_n\PB \mbbh_n^1(\gamma_0) - h_n\PB \mbbh_n^1(\HBN) \\
&= -h_n\int_0^1 \PB^2 \mbbh_n^1(\gamma_u)du \cdot (\HBN - \gamma_0) \\
&=\bigg\{-h_n \PB^2 \mbbh_n^1(\gamma_0) -h_n\int_0^1\int_0^u \PB^3 \mbbh_n^1(\gamma_v) \cdot (\HBN - \gamma_0) dvdu\bigg\} \cdot (\HBN - \gamma_0),
}
where $\gamma_u=u\HBN + (1-u)\gamma_0$.
Together with (\ref{alpha-est-eq}) and the same argument to the proof of Proposition 3.8 in~\cite{ogi23}, we obtain
$$
\begin{aligned}
& \sqrt{T_{n}}\binom{\HBN-\gamma_{0}}{\HAN-\alpha_{0}}
= \MAT{c}{\sqrt{T_n}(T_n\mathcal{I}_\gamma + o_p(T_n))^{-1}h_n\PB \mbbh_n^1(\gamma_0) \\ 
(\mathcal{I}_\alpha)^{-1} T_n^{-1/2}\PA\mbbh_n^2(\alpha_0) + o_p(1) }
=\MAT{c}{(\mathcal{I}_\gamma + o_p(1))^{-1}T_n^{-1/2}h_n\PB \mbbh_n^1(\gamma_0) \\ 
(\mathcal{I}_\alpha)^{-1} T_n^{-1/2}\PA \mbbh_n^2(\alpha_0) + o_p(1) }
\end{aligned}
$$
when the inverse matrix exists.

Then, it is sufficient to show that
\EQ{\label{score-conv-eq} \frac{1}{\sqrt{T_{n}}}\binom{h_{n} \PB \mbbh_n^1(\gamma_0)}{\PA \mbbh_n^2(\alpha_0)} \TOD N\left(0, \Gamma_{1}\right)}
as $n\to\infty$.
(\ref{G1-eq}), (\ref{alscore}) and (\ref{Gn-approx-eq}) yield
$$
\begin{aligned}
T_{n}^{-1/2}\PA \mbbh_n^2(\alpha_0) = & T_{n}^{-1/2} G_{n}^2+o_p(1)=\sum_k\xi_k+o_p(1), \\
h_nT_{n}^{-1/2} \PB \mbbh_n^1(\gamma_0) = & - T_n^{-1/2} G_{n}^1 +o_p(1)=\sum_k\zeta_k+o_p(1).
\end{aligned}
$$
Moreover, we have
$$E[(\zeta_k^\top, \xi_k^\top)^\top (\zeta_k^\top, \xi_k^\top)\mid \FCSK ]
=\frac{1}{T_{n}}\int_{K^k}\int \mathcal{S}\mathcal{S}^\top(t,z) \nu(dz)dt.
$$
Together with (A6), we have
\EQQ{\sum_{k=1}^{\mcl_n} E[(\zeta_k^\top, \xi_k^\top)^\top (\zeta_k^\top, \xi_k^\top) \mid \FCSK ] \TOP \Gamma_1
}
as $n\to\infty$.

Then, (\ref{Gn-approx-eq}), Lemma~\ref{zeta4-conv}, and the martingale central limit theorem (Corollary 3.1 and the remark after that in Hall and Heyde~\cite{hal-hey80}) yield (\ref{score-conv-eq}).

\qed

\subsection{Proof of Proposition~\ref{A6-suff-prop}}

In this subsection, we will show Proposition~\ref{A6-suff-prop}.
We first introduce notation that is needed to decompose each block of $F_n$ into tractable scalar quantities.
With this notation, the proof proceeds as follows.
We decompose $F_n$ into the sums of the terms $\DFN$, $\CFN$ and $\DDFN$, together with boundary terms $R_{n,11}$ and $R_{n,12}$ caused by unequal terminal observation times.
We then localize these quantities on the partition in (B3) and obtain local moment bounds.
After centering the local terms, we apply a mixing inequality under (B4-$1$) to show that the random fluctuations around their expectations vanish.
Finally, (B5) gives the convergence of the expectations, and hence $F_n$ converges in probability to the matrix $\Gamma_1$ required in (A6).

Let $\mathfrak{I}_{l}^{(2)}=\left(\left|I_{i}^{l}\right|\right)_{i=1}^{{\bf M}_{l,n}}$ for $l\in\{1,2\}$,
$$
\begin{aligned}
& \mathcal{G}_{1,1, p, i, i^{\prime}}=[D_{1}^{-1 / 2}\left(G G^{\top}\right)^{p} D_{1}^{-1 / 2}]_{i i}, 
\quad \mathcal{G}_{2,2, p, i, i'}=[D_{2}^{-1 / 2}\left(G^{\top} G\right)^{p} D_{2}^{-1 / 2}]_{i',i'}, \\
& \mathcal{G}_{1,2, p, i, i'}=\mathcal{G}_{2,1, p, i, i'}=[D_{1}^{-1 / 2}\left(G G^{\top}\right)^{p} G D_{2}^{-1 / 2}]_{i, i'}, \\
& S_{1,1, p} = D_1^{-1/2}(GG^\top)^pD_1^{-1/2}, 
\quad S_{2,2, p}=D_2^{-1/2}(G^\top G)^p D_2^{-1/2}, \\
& S_{1,2, p}=\left\{-\mathcal{G}_{1,2, p, i, i^{\prime}}\right\}_{1\leq i\leq \M_{1,n}, 1\leq i^{\prime}\leq \M_{2,n}}, \quad S_{2,1, p}=S_{1,2, p}^{\top},
\end{aligned}
$$
and let
$$
\check{\mathcal{G}}_{m_{1}, m_{2}, q, i, i'}=\DIV{ll}
{\left[S_{m_{1}, m_{2}, q} \mathfrak{I}_{m_{2}}^{(2)}\right]_{i}, &  \text { if } \quad m_1=1, \\
\left[S_{m_{1}, m_{2}, q} \mathfrak{I}_{m_{2}}^{(2)}\right]_{i^{\prime}},  & \text { if } \quad m_1=2}
$$
for nonnegative integer $p$ and $q$.

Let
\EQNN{
\check{F}_{l_1,l_2,m_1,m_2, p,q,n}&=\frac{h_{n}}{T_{n}} \SUMIIp |I_i\cap I_{i'}|\mathcal{G}_{l_1,l_2,p,i,\varphi(i')}\check{\mathcal{G}}_{m_1,m_2,q,i,\varphi(i')}, \\
\ddot{F}_{l_1,l_2,p,n}&=\frac{1}{T_{n}}(\mathfrak{I}_{l_{1}}^{(2)})^{\top} S_{l_{1}, l_{2}, p} \mathfrak{I}_{l_{2}}^{(2)}, \\
\dot{F}_{l_{1}, l_{2}, m_{1}, m_{2}, p, q,n}&=\frac{h_{n}^2}{T_{n}} \SUMIIp
|I_{i} \cap I_{i'}| \mcg_{l_{1}, l_{2}, p, i, \varphi(i')} \mcg_{m_1, m_{2}, q, i, \varphi(i')}.
}

Let $e_{a,b}=1$ if $a=b$ and $e_{a,b}=-1$ if $a\ne b$, and set $C_0=C(\gamma_0)$,
$$a_{l, l, p}(\gamma)=\rho^{2 p}(V^l)^{-1}(\gamma), \quad a_{1,2, p}(\gamma)=a_{2,1, p}(\gamma)=\rho^{2 p+1}(V^1V^2)^{-1/2}(\gamma),$$
$$\mathcal{A}_{l_1, l_2, m} = \frac12 \int[C_0z]_{l_1}[C_0z]_{l_2}[C_0z]_{m} \nu(d z),$$
\EQN{\mca_{l_1,l_2,m_1,m_2} = \frac14 \int [C_0z]_{l_1}[C_0z]_{l_2}[C_0z]_{m_1}[C_0z]_{m_2} \nu(dz).}

Let
$$
\begin{aligned}
\dot{\mca}_{l_1, l_{2}, m_{1}, m_{2}, p, q} 
= \left(e_{l_1,l_2}e_{m_1,m_2} \mca_{l_{1}, l_{2}, m_{1}, m_{2}} \partial_{\gamma_i} a_{l_{1}, l_{2}, p}(\gamma_0) \partial_{\gamma_j} a_{m_1, m_{2}, q}(\gamma_0) \right)_{1\leq i,j\leq p_\gamma},
\end{aligned}
$$
$$
\begin{aligned}
& \check{\mathcal{A}}_{l_{1}, l_{2}, m_{1}, m_{2}, p,q}  = \left( - e_{l_1,l_2}\mca_{l_{1}, l_{2}, m_1} a_{m_{1}, m_{2}, q}(\gamma_0) \partial_{\alpha_j} A^{m_{2}}(\alpha_0) \partial_{\gamma_i} a_{l_{1}, l_{2}, p}(\gamma_0) \right)_{1\leq i\leq p_\gamma, 1\leq j\leq p_\alpha},
\end{aligned}
$$
$$
\ddot{\mca}_{l_1,l_2,p} = \left(a_{l_{1}, l_{2}, p}(\gamma_0) \partial_{\alpha_i} A^{l_{1}}(\alpha_0)  \partial_{\alpha_j} A^{l_{2}}(\alpha_0) \right)_{1\leq i,j\leq p_\alpha}.
$$

Let $\Gamma_1$ denote the matrix in (B5). We decompose $F_n$ as 
$$F_n=\MAT{cc}{F_{n,11} & F_{n,12} \\ F_{n,12}^\top & F_{n,22}}.$$

To describe the boundary terms, let
$$\tau_{l,n}=S_{{\bf M}_{l,n}}^{n,l},\qquad
d_{i,n}=|I_i\cap[\tau_{3-\psi(i),n},T_n)|
\quad (1\leq i\leq {\bf M}_n).$$
Thus, $d_{i,n}$ is the length of the part of $I_i$ not covered by observation intervals of the other component.
Since the intervals of component $l$ partition $[0,\tau_{l,n})$, we have
\EQ{\label{terminal-overlap-eq}
|I_i|=\sum_{\substack{1\leq j\leq {\bf M}_n\\\psi(j)\ne\psi(i)}}|I_i\cap I_j|+d_{i,n},
\qquad \sum_{i=1}^{{\bf M}_n}d_{i,n}=|\tau_{1,n}-\tau_{2,n}|.
}
Let $\mathcal C_k=\partial_{\gamma_k}S_n^{-1}(\gamma_0)$ and
$\mathcal B_a=S_n^{-1}(\gamma_0)\partial_{\alpha_a}\Delta\mathcal V(\alpha_0)$.
Define the $p_\gamma\times p_\gamma$ matrix $R_{n,11}$ and the $p_\gamma\times p_\alpha$ matrix $R_{n,12}$ by
\EQN{\label{Fn-boundary-def}
[R_{n,11}]_{k,k'}&=\frac{h_n^2}{T_n}\sum_{i=1}^{{\bf M}_n}
d_{i,n}\mathcal A_{\psi(i),\psi(i),\psi(i),\psi(i)}
[\mathcal C_k]_{ii}[\mathcal C_{k'}]_{ii},\\
[R_{n,12}]_{k,a}&=-\frac{h_n}{T_n}\sum_{i=1}^{{\bf M}_n}
d_{i,n}\mathcal A_{\psi(i),\psi(i),\psi(i)}
[\mathcal C_k]_{ii}[\mathcal B_a]_i.
}

Both boundary terms vanish when $\tau_{1,n}=\tau_{2,n}$.

\begin{lemma}\label{Fn-lem}
Assume (A1) and (A2). Then,
\EQNN{ F_{n,11} &= \SUMPQ \SUMLL \SUMMM
\dot{\mca}_{l_1,l_2,m_1,m_2,p,q} \dot{F}_{l_1,l_2,m_1,m_2,p,q,n}+R_{n,11}, \\
F_{n,12} &= \SUMPQ \SUMLL \SUMMM
\check{\mathcal{A}}_{l_1,l_2,m_1,m_2, p,q} \check{F}_{l_1,l_2,m_1,m_2,p,q,n}+R_{n,12}, \\
F_{n,22} &= \SUMP \SUMLL \ddot{\mca}_{l_1,l_2,p} \ddot{F}_{l_1,l_2,p,n}.
}
\end{lemma}

\begin{proof}

We first expand $F_{n,11}$ without replacing the lengths of individual intervals by sums of overlaps.

\EQNN{[F_{n,11}]_{k,l} &= \frac{h_n^2}{4T_n}
 \int_0^{T_n}\int \chi_{Cz}(t)^{\top} \mcc_k\chi_{Cz}(t) \chi_{Cz}(t)^\top \mcc_l \chi_{Cz}(t) \nu(dz) d t \\
&= \frac{h_n^2}{4T_n} \sum_{i_1, \cdots i_4} \left[\mcc_k\right]_{i_1, i_2} \left[\mcc_l\right]_{i_3,i_4} \int_0^{T_n} \int  \prod_{j=1}^{4} [\chi_{Cz}(t)]_{i_j} \nu (dz)dt.
}

$\int_0^{T_n} \prod_{j=1}^{4} [\chi_{Cz}(t)]_{i_j} d t$ is non-zero only if $i_1=\cdots = i_4$ or the case that $\{i_1,\cdots,i_4\}$ is equal to a set of indicators $\{i,i'\}$ with $i\leq {\bf M}_{1,n}<i'$ and $I_i\cap I_{i'}\neq \emptyset$.
Therefore, by setting $\mfa^{(k)}_{i,i'}(z)=[Cz]_{\psi(i)}[Cz]_{\psi(i')}[\mcc_k]_{i,i'}$, we obtain
\EQN{\label{zeta2-eq}
[F_{n,11}]_{k,l} &= \frac{h_n^2}{4T_{n}}\sum_{i}|I_i| \int \mfa^{(k)}_{i,i} \mfa^{(l)}_{i,i}\nu(dz) \\
&\quad +\frac{h_n^2}{4T_{n}} \SUMIIp|I_i\cap I_{i'}|
\int \Big\{ 2\mfa^{(k)}_{i,i'}(\mfa^{(l)}_{i,i} + \mfa^{(l)}_{i',i'}) + 2\mfa^{(l)}_{i,i'}(\mfa^{(k)}_{i,i} + \mfa^{(k)}_{i',i'}) \\
&\qquad \qquad \qquad \qquad \qquad + 4\mfa^{(k)}_{i,i'}\mfa^{(l)}_{i,i'} + \mfa^{(k)}_{i,i}\mfa^{(l)}_{i',i'} + \mfa^{(l)}_{i,i}\mfa^{(k)}_{i',i'} \Big\}\nu(dz). 
}

Here the coefficients in the second term of (\ref{zeta2-eq}) come from the possible positions of the two indices corresponding to the off-diagonal block.
The products of the type $(i,i')$ and $(i,i)$ give four possible arrangements, which are summarized as the two terms
$2\mfa^{(k)}_{i,i'}\mfa^{(l)}_{i,i}$ and $2\mfa^{(l)}_{i,i'}\mfa^{(k)}_{i,i}$.
Similarly, the products of the type $(i,i')$ and $(i',i')$ give the two terms.
The products of two off-diagonal terms $(i,i')$ and $(i,i')$ give $2\times 2$ arrangements, yielding the coefficient $4$.
Finally, the products of the two diagonal terms $(i,i)$ and $(i',i')$ appear in two orders, which gives the last two terms in the braces.

\begin{note}
{\colorr 
$(i,i'),(i,i)$: 4 cases depending on location of $i'$,\\
$(i,i'),(i',i')$: 4 cases depending on location of $i$,\\
$(i,i'),(i,i')$: $2\times 2$ cases depending on location of $i$,
\\$(i,i),(i',i')$: 2 cases depending on location of $i$.
}
\end{note}

Lemma~\ref{ogi23-lem} yields $\|\tilde{G}\|<1$.
Then, (\ref{invSn-eq}) yields
\EQN{\label{invSn-eq2}
S_{n}^{-1}(\gamma)&=\tilde{D}^{-1/2} \sum_{p=0}^{\infty} \MAT{cc}{(\tilde{G}\tilde{G}^\top)^p & -(\tilde{G}\tilde{G}^\top)^p\tilde{G} \\
-(\tilde{G}^\top\tilde{G})^p\tilde{G}^{\top} & (\tilde{G}^\top\tilde{G})^p}
\tilde{D}^{-1/2} \\
&=\sum_{p=0}^{\infty}\MAT{cc}{
a_{1,1, p} S_{1,1, p} & a_{1,2, p} S_{1,2, p} \\
a_{2,1, p} S_{2,1, p} & a_{2,2, p} S_{2,2, p}}.
}

For each fixed $n$, the series in (\ref{invSn-eq2}) and its first
$\gamma$-derivatives converge absolutely and uniformly in
$\gamma\in\overline{\Theta}_\gamma$, almost surely.
Indeed, (A2) implies that the coefficients with $p=0$ and their first
derivatives are bounded, and that
\[
 \sup_{\gamma\in\overline{\Theta}_\gamma}
 \left\{|a_{l_1,l_2,p}(\gamma)|
       +|\partial_\gamma a_{l_1,l_2,p}(\gamma)|\right\}
 \leq C(p+1)\bar\rho^{\,2p-1},\qquad p\geq1,
\]
where $C$ is independent of $p$ and $\bar\rho<1$.
Since $\|G\|\leq1$ by Lemma~\ref{ogi23-lem}, we also have
$\|S_{l_1,l_2,p}\|\leq\underline r_n^{-1}<\infty$ for every $p$.
The assertion therefore follows from
$\sum_{p=1}^{\infty}(p+1)\bar\rho^{\,2p-1}<\infty$.
Consequently, (\ref{invSn-eq2}) can be differentiated term by term, and products of
the resulting series are absolutely convergent double series.
These bounds justify the rearrangements and the exchanges of the series
with the integrals used below in the decompositions of
$F_{n,11}$, $F_{n,12}$ and $F_{n,22}$.

Using (\ref{terminal-overlap-eq}) for the terms with all four component indices equal,
together with (\ref{zeta2-eq}), we obtain
\EQN{& [F_{n,11}]_{k,l} \\
&\quad =\frac{h_n^2}{T_{n}}\sum_{i\leq {\bf M}_{1,n}}\mca_{1,1,1,1}|I_i|[\mcc_k]_{i,i}[\mcc_l]_{i,i}
+ \frac{h_n^2}{T_{n}}\sum_{i'> {\bf M}_{1,n}} \mca_{2,2,2,2}|I_{i'}|[\mcc_k]_{i',i'}[\mcc_l]_{i',i'} \\
&\qquad +\frac{h_n^2}{T_{n}} \SUMIIp|I_i\cap I_{i'}|\big\{ 2\mca_{1,1,1,2} ([\mcc_k]_{i,i}[\mcc_l]_{i,i'} + [\mcc_l]_{i,i}[\mcc_k]_{i,i'})
+2\mca_{1,2,2,2}([\mcc_k]_{i',i'}[\mcc_l]_{i,i'} + [\mcc_l]_{i',i'}[\mcc_k]_{i,i'}) \\
&\qquad \qquad \qquad\qquad 
+ 4\mca_{1,1,2,2}[\mcc_k]_{i,i'}[\mcc_l]_{i,i'} 
+ \mca_{1,1,2,2}([\mcc_k]_{i,i}[\mcc_l]_{i',i'} + [\mcc_l]_{i,i}[\mcc_k]_{i',i'}) \big\} \\
 &\quad = \frac{h_{n}^2}{T_{n}} \SUMIIp |I_i\cap I_{i'}| \sum_{p,q=0}^\infty \SUMLL \SUMMM
 \mca_{l_1,l_2,m_1,m_2} \\
 &\qquad\qquad\times e_{l_1,l_2}e_{m_1,m_2}\partial_{\gamma_k} a_{l_1,l_2,p} \partial_{\gamma_l} a_{m_1,m_2,q} 
 \mcg_{l_1,l_2,p,i,\varphi(i')} \mcg_{m_1,m_2,q,i,\varphi(i')}+[R_{n,11}]_{k,l} \\
&\quad =\SUMPQ \SUMLL \SUMMM [\dot{\mca}_{l_1, l_2, m_{1}, m_{2}, p, q}]_{k,l} \dot{F}_{l_1, l_2, m_{1}, m_{2}, p, q,n}+[R_{n,11}]_{k,l}.
}

For the cross block, using $\mathcal B_l$ defined above, we obtain
$$
\begin{aligned}
{}[F_{n,12}]_{k,l}
&=-\frac{h_{n}}{2T_{n}} \int_0^{T_n}\int
\chi_{Cz}(t)^\top\mcc_k\chi_{Cz}(t)\mcb_l^\top\chi_{Cz}(t)\nu(dz)dt \\
&=-\frac{h_{n}}{T_{n}} \sum_{i \leq {\bf M}_{1, n}} \mathcal{A}_{1,1,1}\left|I_{i}\right|\left[\mcc_k\right]_{i i}\left[\mcb_l\right]_{i} 
-\frac{h_{n}}{T_{n}} \sum_{i^{\prime}>{\bf M}_{1, n}} \mathcal{A}_{2,2,2}\left|I_{i^{\prime}}\right|\left[\mcc_k\right]_{i' i'}\left[\mcb_l\right]_{i^{\prime}} \\
 & \quad - \frac{h_{n}}{T_{n}}\SUMIIp|I_i\cap I_{i'}|
 \bigg\{\mathcal{A}_{1,2,2}\left[\mcc_k\right]_{i' i'}\left[\mcb_l\right]_{i}
 +\mathcal{A}_{1,1,2}\left[\mcc_k\right]_{i i}\left[\mcb_l\right]_{i'} \\
&\qquad \qquad +2 \mathcal{A}_{1,1,2}\left[\mcc_k\right]_{i i'}\left[\mcb_l\right]_{i}
+2 \mathcal{A}_{1,2,2}\left[\mcc_k\right]_{i i'}\left[\mcb_l\right]_{i'} \bigg\} \\
&=-\frac{h_{n}}{T_{n}}\sum_{\substack{i\leq {\bf M}_{1,n} \\ i'>{\bf M}_{1,n}}}|I_i\cap I_{i'}|
\SUMLL \SUMM \mathcal{A}_{l_1,l_2,m}\mathfrak{A}^{(k)}_{l_1,l_2,i,i'}\tilde{\mathfrak{A}}^{(l)}_{m,i,i'}+[R_{n,12}]_{k,l},
\end{aligned}
$$
where
$$\mathfrak{A}^{(k)}_{1,1, i, i^{\prime}}=\left[\mcc_k\right]_{i i}, 
\mathfrak{A}^{(k)}_{1,2, i, i'}=\mathfrak{A}^{(k)}_{2,1, i, i'}=\left[\mcc_k\right]_{i i'}, 
\mathfrak{A}^{(k)}_{2,2, i, i'}=\left[\mcc_k\right]_{i',i'},$$
$$
\tilde{\mathfrak{A}}^{(l)}_{1, i, i'}=\left[\mcb_l\right]_{i}, 
\quad \tilde{\mathfrak{A}}^{(l)}_{2, i, i'}=\left[\mcb_l\right]_{i'}.
$$
In the last equality we again used (\ref{terminal-overlap-eq}).
Then, together with (\ref{invSn-eq2}), we obtain
\EQQ{[F_{n,12}]_{k,l} = -\frac{h_n}{T_n}\SUMIIp |I_i\cap I_{i'}|\SUMP\SUMLL\SUMM e_{l_1,l_2}\mca_{l_1,l_2,m}\tilde{\mathfrak{A}}^{(l)}_{m,i,i'}\partial_{\gamma_k}a_{l_1,l_2,p}(\gamma_0) \mathcal{G}_{l_1,l_2,p,i,i'}+[R_{n,12}]_{k,l}.
}

Moreover, since
$\PA \Delta \mcv =(\PA A^1 (\mathfrak{I}_{1}^{(2)})^\top, \PA A^{2} (\mathfrak{I}_{2}^{(2)})^\top )^\top$,
we have
$$ S_{n}^{-1} \PA \Delta \mcv=\sum_{q=0}^{\infty}\binom{a_{1,1, q} \partial_\alpha A^{1} S_{1,1, q} \mathfrak{I}_{1}^{(2)}+a_{1,2, q} \partial_\alpha A^{2} S_{1,2, q} \mathfrak{I}_{2}^{(2)}}{a_{2,1, q} \partial_\alpha A^1 S_{2,1, q} \mathfrak{I}_{1}^{(2)}+a_{2, 2, q} \partial_\alpha A^{2} S_{2,2, q} \mathfrak{I}_{2}^{(2)}}. $$
Then, we have $\tilde{\mathfrak{A}}^{(l)}_{m_1,i,i'}=\sum_{q=0}^{\infty} \sum_{m_{2}=1}^{2} a_{m_{1}, m_{2}, q} \partial_{\alpha_l} A^{m_{2}} \check{\mathcal{G}}_{m_1, m_{2}, q, i, \varphi(i')}$, and hence we obtain
$$
\begin{aligned}
[F_{n,12}]_{k,l}
&=-\frac{h_{n}}{T_{n}} \SUMIIp\left|I_{i} \cap I_{i^{\prime}}\right| \SUMPQ \SUMLL \SUMMM
e_{l_1,l_2}\mca_{l_1, l_2, m_1} a_{m_{1}, m_{2}, q}(\gamma_0)
\partial_{\alpha_l} A^{m_{2}}(\alpha_0) \partial_{\gamma_k} a_{l_1,l_2,p}(\gamma_0) \\
&\qquad \qquad \qquad \qquad \times \mathcal{G}_{l_{1}, l_{2}, p, i, \varphi(i')} \check{\mathcal{G}}_{m_{1}, m_{2}, q, i, \varphi(i')}+[R_{n,12}]_{k,l} \\
&= \SUMPQ \SUMLL \SUMMM
[\check{\mathcal{A}}_{l_1,l_2,m_1,m_2,p,q}]_{k,l}\check{F}_{l_1,l_2,m_1,m_2,p,q,n}+[R_{n,12}]_{k,l}.
\end{aligned}
$$

Similarly, we obtain
$$
\begin{aligned}
[F_{n,22}]_{k,l}
& =\frac{1}{T_{n}} \int_{0}^{T_{n}} \int
\left\{\left(\partial_{\alpha_k}  \Delta \mcv\right)^\top S_{n}^{-1} \chi_{C z}(t)\right\}
\left\{\left(\partial_{\alpha_l}  \Delta \mcv\right)^\top S_{n}^{-1} \chi_{C z}(t)\right\} d t \nu(d z) \\
& =\frac{1}{T_n}\left(\partial_{\alpha_k} \Delta \mcv\right)^{\top} S_{n}^{-1} \partial_{\alpha_l} \Delta \mcv \\
& =\frac{1}{T_{n}} \SUMP
\left(\partial_{\alpha_k} A^{1}(\mathfrak{I}_{1}^{(2)})^{\top}, \partial_{\alpha_k} A^{2}(\mathfrak{I}_{2}^{(2)})^{\top}\right)
\left(\begin{array}{ll}
a_{1,1, p} S_{1,1, p} & a_{1,2, p} S_{1,2, p} \\
a_{2,1, p} S_{2,1, p} & a_{2,2,p} S_{2,2, p}
\end{array}\right) 
 \binom{\partial_{\alpha_l} A^{1} \mathfrak{I}_{1}^{(2)}}{\partial_{\alpha_l} A^{2} \mathfrak{I}_{2}^{(2)}} \\
& = \SUMP \SUMLL [\ddot{\mca}_{l_1,l_2,p}]_{k,l}\ddot{F}_{l_1,l_2,p,n}.
\end{aligned}
$$

\end{proof}

The proof of the following lemma is left to the appendix.
\begin{lemma}\label{GGTp-lem}
For any $p\in \mbbz_+$, $1\leq i\leq \M_{1,n}$ and $1\leq i'\leq \M_{2,n}$, 
    $$[(GG^\top)^pG]_{ii'}\frac{1}{\sqrt{|I_{i'}|}}\leq \frac{1}{\sqrt{|I_{i}|}},
    \quad [(G^\top G)^pG^\top]_{i'i}\frac{1}{\sqrt{|I_i|}}\leq \frac{1}{\sqrt{|I_{i'}|}}.$$
\end{lemma}

For $x\in \mathbb{R}$, let $[x]$ denote the greatest integer not exceeding $x$.
For $(t_k^n)_{k=1}^{L_n}$ in (B3), 
let $\tau_{k}=t_{k}^n-t_{k-1}^n$, and let $\tau_{k, p}=\tau_k+(2 p+2) r_{n}$, $\underline{\tau}=\min_{1\leq k\leq L_n}\tau_k$ and $\bar{\tau}=\max_{1\leq k\leq L_n}\tau_k$.
By (B3), there exists a positive constant $C_\tau$ such that
$\bar{\tau}/\underline{\tau}\leq C_\tau$ for all $n$.
Since $T_n=\sum_{k=1}^{L_n}\tau_k$, we have
$$L_nT_n^{-1}\underline{\tau}
\geq \frac{L_n\underline{\tau}}{L_n\bar{\tau}}\geq C_\tau^{-1},
\quad
L_nT_n^{-1}\bar{\tau}
\leq \frac{L_n\bar{\tau}}{L_n\underline{\tau}}\leq C_\tau.
$$
In particular, there exist positive constants $c$ and $C$ such that
\EQ{\label{LnTn-ineq}
0<c\leq L_nT_n^{-1}\underline{\tau}
\leq L_nT_n^{-1}\bar{\tau}\leq C
}
for any $n$.
Then, for $\epsilon\in (0,1)$ in (B3), we have
\EQ{\label{LnTn-ineq2} 
n^{-\epsilon}\underline{\tau}h_n^{-1}
=L_nT_n^{-1}\underline{\tau}\cdot \frac{n^{1-\epsilon}}{L_n}\to \infty.}

Let
$$
\begin{aligned}
& \mci_k^{l}=\left\{1 \leq i \leq \M_{l,n}; I_{i}^{l} \cap\left[t_{k-1}^n, t_k^n\right) \neq \phi\right\}, \quad \mci_{k}=\mci_{k}^{1} \cup \{i+\M_{1,n}; i\in \mci_k^2\}.
\end{aligned}
$$

Let $\mathfrak{I}_{l}^{(2,k)}=(|I_i^l\cap [t_{k-1}^n,t_k^n)|)_{i=1}^{\M_{l,n}}$, and let
\EQNN{\CFK&=\frac{1}{n} \SUMIIp |I_i\cap I_{i'}\cap \TKI|\mathcal{G}_{l_1,l_2,p,i,\varphi(i')}\check{\mathcal{G}}_{m_1,m_2,q,i,\varphi(i')}, \\
\DDFK&=\frac{1}{T_{n}}(\mathfrak{I}_{l_{1}}^{(2)})^{\top} S_{l_{1}, l_{2}, p} \mathfrak{I}_{l_{2}}^{(2,k)}, \\
\DFK&=\frac{h_{n}}{n} \SUMIIp
|I_{i} \cap I_{i'}\cap \TKI| \mcg_{l_{1}, l_{2}, p, i, \varphi(i')} \mcg_{m_1, m_{2}, q, i, \varphi(i')}.}

Then, we obtain the following bounds for $\ddot{F}_{l_1,l_2,p,n,k}$, $\check{F}_{l_1,l_2,m_1,m_2, p,q,n,k}$ and
$\dot{F}_{l_1,l_2,m_1,m_2,p,q,n,k}$.
The proof is left to the appendix.
\begin{lemma}\label{three-F-lem}
Assume (A1), (A2) and (B3). Then, there exists a positive constant $C$ such that 
for any $l_1,l_2,m_1,m_2 \in \{1,2\}$, $p,q\in \mbbz_+$, $1\leq k\leq L_n$ and $n\in \mbbn$,
\EQ{\label{mca-eq} \|\dot{\mca}_{l_1, l_{2}, m_{1}, m_{2}, p, q}\|\vee \|\check{\mathcal{A}}_{l_{1}, l_{2}, m_{1}, m_{2}, p,q}\|
 \leq C\bar{\rho}^{(2p+2q-2)\vee 0},
\quad \|\ddot{\mca}_{l_1,l_2,q}\| \leq C\bar{\rho}^{2q},
}
$$|\ddot{F}_{l_1,l_2,q,n,k}| \leq T_n^{-1}\tau_{k,q},$$
\EQQ{|\check{F}_{l_1,l_2,m_1,m_2, p,q,n,k}| \leq \frac{\sqrt{\tau_{k,q}}}{n}\sqrt{\sum_{i\in \mci_k}|I_i|^{-1}}, 
\quad |\dot{F}_{l_1,l_2,m_1,m_2,p,q,n,k}| \leq \frac{h_n}{n}
\sum_{i\in \mci_k}|I_i|^{-1}.
}
Moreover, there exists a positive constant $C'$ such that 
$$E[\ddot{F}_{l_1,l_2,q,n,k}^4]\vee E[\check{F}_{l_1,l_2,m_1,m_2, p,q,n,k}^4]\vee E[\dot{F}_{l_1,l_2,m_1,m_2,p,q,n,k}^4]
\leq C' T_n^{-4}(\tau_k^4+(q^4+1)E[r_n^4]) $$
for any $l_1,l_2,m_1,m_2,p,q,k$ and $n$.
\end{lemma}

\noindent
{\bf Proof of Proposition~\ref{A6-suff-prop}.}

Let $\bar{E}[{\bf X}]={\bf X} - E[{\bf X}]$ for a random variable ${\bf X}$, and let $R_{p,q}=2(p\vee q)+3$.
For $\epsilon\in (0,1)$ in (B3), let
$$
\begin{aligned}
& A_{h, t}^{p,+}=\cap_{i=1,2} \cap_{l \in\left[1, p \wedge h^{-1}(T_n-t)\right] \cap \mathbb{N}}\left\{\omega ; \mcn_{t+l h}^i-\mcn_{t+(l-1) h}^i>0\right\}, \\
& A_{h, t}^{p,-}=\cap_{i=1,2} \cap_{l \in\left[1, p \wedge h^{-1} t\right] \cap \mathbb{N}}\left\{\omega ; \mcn_{t-(l-1) h}^i-\mcn_{t-l h}^i>0\right\}, \\
& \tilde{A}_{k, \epsilon}^{p,q}=A_{h_n[n^\epsilon], t_k^n}^{R_{p,q},+} \cap A_{h_n[n^\epsilon], t_{k-1}^n}^{R_{p,q},-}.
\end{aligned}
$$
Let $X'_k$ be one of $\bar{E}[\check{F}_{l_1,l_2,m_1,m_2, p,q,n,k} 1_{\tilde{A}_{k,\epsilon}^{p,q}}]$,
$\bar{E}[\ddot{F}_{l_1,l_2,p,n,k} 1_{\tilde{A}_{k,\epsilon}^{p,q}}]$,
or $\bar{E}[\dot{F}_{l_1,l_2,m_1,m_2,p,q,n,k} 1_{\tilde{A}_{k,\epsilon}^{p,q}}]$.

Let $\delta_k^-=\max\{i\in \mbbz_+ ; ih_n\leq t_{k-1}^n\}$
and $\delta_k^+=\min\{i\in \mbbz_+ ; ih_n\geq t_k^n\}$.
If \(\omega \in \tilde{A}_{k,\epsilon}^{p,q}\), for the sequence of intervals 
$$\{(t_k^n + (l-1)h_n [n^{\epsilon}], t_k^n + l h_n [n^{\epsilon}]]\}_{l=1}^{R_{p,q}} \quad  {\rm and} \quad  \{(t_{k-1}^n - l h_n [n^{\epsilon}], t_{k-1}^n - (l-1) h_n [n^{\epsilon}]]\}_{l=1}^{R_{p,q}}, 
$$
each interval includes at least one jump point of both $\mcn_t^1(\omega)$ and $\mathcal{N}_t^2(\omega)$ 
if the intervals are contained within \([0, T_n]\).

On $\widetilde A_{k,\epsilon}^{p,q}$, each overlap path
contributing to the three local quantities has at most
$2(p\vee q)+1$ edges. Since each full coarse interval
contains an observation time of both components, each step
along such a path extends the required window by at most
one coarse interval. Thus, with $R_{p,q}=2(p\vee q)+3$,
all observation intervals needed to compute these quantities
can be reconstructed from the counting-process increments
in the enlarged window. At the boundaries, we use the known
initial time $0$ and only the actual observation times up to
$T_n$.
Since the local quantities are invariant under relabelling
of the observation intervals, and
$\widetilde A_{k,\epsilon}^{p,q}$ is measurable with respect
to the same increments, their products with
$1_{\widetilde A_{k,\epsilon}^{p,q}}$, and hence $X'_k$,
are $\mathcal H^n_{t(k,p,q,-),t(k,p,q,+)}$-measurable,
where
$$t(k,p,q,+) = (\delta_k^++R_{p,q}[n^{\epsilon}])\wedge n,
\quad t(k,p,q,-) = (\delta_k^--R_{p,q}[n^{\epsilon}])\vee 0.
$$

Let $\alpha'_0=1/4$, $\alpha'_k=0$ for $k\geq L_n$ and 
$$\alpha'_k=\sup_{\substack{1\leq i<j\leq L_n \\j-i\geq k}} 
\sup_{\substack{A\in \sigma(X'_l;l\leq i) \\ B\in \sigma(X'_m;m\geq j)}}|P(A\cap B)-P(A)P(B)|$$
for $1\leq k<L_n$, and let $t(p,q)=R_{p,q}[n^{\epsilon}]$.
Then, since 
$$(t(j,p,q,-)-t(i,p,q,+))\vee 0\geq \delta_j^--\delta_i^+-2t(p,q)
\geq (t_j^n-t_i^n)h_n^{-1}-2-2t(p,q)\geq [(k-1)\underline{\tau}h_n^{-1}]-2-2t(p,q)$$
for $2\leq k<L_n$ and $1\leq i\leq i+k\leq j\leq L_n$, we have $\alpha'_k\leq \alpha_{([(k-1)\underline{\tau}h_n^{-1}]-2t(p,q)-2)\vee 0}^n$.
\begin{note}
{\colorr If $t(j,p,q,-)=0$ or $t(i,p,q,+)=n$,
then the inequality obviously satisfied.
If $t(j,p,q,-)\leq t(i,p,q,+)$, the RHS is $\alpha_0^n=1/4$. So the inequality is always satisfied.} 
\end{note}

Let $K=\inf_{k\geq 1}([k\underline{\tau}h_n^{-1}]-[(k-1)\underline{\tau}h_n^{-1}])$.
Then, when $u<1/4$, Lemma~\ref{alpha-k-lem} in the appendix yields
\EQQ{\alpha'^{-1}(u)=\sum_{k=0}^\infty 1_{\{\alpha'_k>u\}}
\leq 2+K^{-1}(\alpha^{-1}(u)+2t(p,q)+1).
}

Then, the Rosenthal-type inequality (Theorem 4 in Ogihara and Yoshida \cite{ogi-yos14}) yields that
\EQNN{E\bigg[\bigg|\sum_{k=1}^{L_n}X'_k\bigg|^2\bigg]
&\leq C\sum_{k=1}^{L_n}\int_0^{1/4}(2+K^{-1}(\alpha^{-1}(u)+2t(p,q)+1))Q_{X'_k}^2(u)du \\
&\leq C\sum_{k=1}^{L_n}\sqrt{\int_0^1(2+K^{-1}(\alpha^{-1}(u)+2t(p,q)+1))^2du} \sqrt{\int_0^1Q_{X'_k}^4(u)du},
}
where $Q_{X'_k}(u)=\inf\{t>0; P(|X'_k|>t)\leq u\}$.

By a similar argument to the proof of Proposition 6 in \cite{ogi-yos14} and (\ref{alpha-k-eq}), we obtain
$$\int_0^1Q_{X'_k}^4(u)du=E[|X'_k|^4], 
\quad \int_0^1\left(\alpha^{-1}(u)\right)^l du\leq 2\sum_{k=0}^\infty (k+1)^{l-1}\alpha_k^n<\infty$$
for $l=1,2$.
Together with Lemma~\ref{three-F-lem}, we obtain
\EQNN{E\bigg[\bigg|\sum_{k=1}^{L_n}X'_k\bigg|^2\bigg]
&\leq CL_n\sqrt{1+K^{-2}t(p,q)^2}\sup_k\sqrt{E[|X'_k|^4]} \\
&\leq C_{p,q}L_nT_n^{-2}(\bar{\tau}^2+E[r_n^4]^{1/2})\sqrt{1+K^{-2}t(p,q)^2}.
}
for some positive constants $C$ and $C_{p,q}$.

Since (\ref{LnTn-ineq2}) yields
$$K^{-2}t(p,q)^2
\leq C'_{p,q}h_n^2\underline{\tau}^{-2}n^{2\epsilon}\to 0$$
as $n\to\infty$ for some positive constant $C'_{p,q}$.
Together with (B3) and (\ref{LnTn-ineq}), we have
\EQ{\label{tXk-conv} E\bigg[\bigg|\sum_{k=1}^{L_n}X'_k\bigg|^2\bigg]\leq C(L_n^2T_n^{-2}\bar{\tau}^2+o(1))\cdot L_n^{-1}\to 0.}

(B2-$r$) implies
$$
P((\tilde{A}_{k,\epsilon}^{p,q})^{c}) \leq (4 p \vee q+2) \max _{l=1,2} P\left(\mcn_{t+h_{n}\left[n^{\epsilon}\right]}^{n, l}-\mcn_{t}^{n, l}=0\right) 
\leq c_{p, q} n^{-\epsilon r}.
$$
Then, Lemma~\ref{three-F-lem}, (\ref{LnTn-ineq}) and (B3) yield
\EQNN{E\bigg[\sum_{k=1}^{L_n}\DFK 1_{(\TAK)^c}\bigg]
&\leq \sum_{k=1}^{L_n}E[\DFK^2]^{1/2}P((\TAK)^c)^{1/2} \\
&\leq CL_nT_n^{-1}(\bar{\tau}+(q+1)E[r_n^4]^{1/4})n^{-\epsilon r/2}\to 0.
}
Together with (\ref{tXk-conv}), we obtain
$\bar{E}[\DFN]=\bar{E}[\sum_k\DFK]\TOP 0$.
Similarly, we obtain
$\bar{E}[\DDFN] \TOP 0$ and $\bar{E}[\CFN]\TOP 0$.

Moreover, for any $\delta>0$, Lemma~\ref{three-F-lem} and (\ref{LnTn-ineq}) yield
\EQNN{&E\bigg[\bigg\|\sum_{p,q;p\geq M~{\rm or}~q\geq M}\dot{\mca}_{l_1,l_2,m_1,m_2,p,q}\DFN \bigg\| \bigg] \\
&\quad \leq C\sum_{p,q;p\geq M~{\rm or}~q\geq M}\bar{\rho}^{(2p+2q-2)\vee 0}\sum_{k=1}^{L_n}T_n^{-1}(\tau_k+(q+1) E[r_n^4]^{1/4}) \\
&\quad \leq C\sum_{p,q;p\geq M~{\rm or}~q\geq M}\bar{\rho}^{(2p+2q-2)\vee 0}(1+(q+1) L_nT_n^{-1}E[r_n^4]^{1/4}) < \delta
}
for sufficiently large $M$.

It remains to control the boundary terms in Lemma~\ref{Fn-lem}.
Lemmas~\ref{ogi23-lem} and \ref{invSdelV-est-lemma} give
$|[\mathcal C_k]_{ii}|\leq C|I_i|^{-1}$ and $|[\mathcal B_a]_i|\leq C$.
Set
$$D_n^{\partial}=\frac{h_n^2}{T_n}\sum_{i=1}^{{\bf M}_n}\frac{d_{i,n}}{|I_i|^2},
\qquad Y_{n,k}^{\partial}=\frac{h_n^2}{T_n}\sum_{i\in\mathcal I_k}|I_i|^{-1}.$$
By (\ref{Fn-boundary-def}), (\ref{terminal-overlap-eq}) and the Cauchy--Schwarz inequality,
$$\|R_{n,11}\|\leq CD_n^{\partial},\qquad
\|R_{n,12}\|\leq C(D_n^{\partial})^{1/2}
\left(\frac{|\tau_{1,n}-\tau_{2,n}|}{T_n}\right)^{1/2}
\leq C(D_n^{\partial})^{1/2}.$$
Let $E_n^{\partial}=\{\min_{l=1,2}\tau_{l,n}\geq T_n-\underline\tau\}$.
On this event, every interval with $d_{i,n}>0$ intersects the last block $[t_{L_n-1}^n,T_n)$, so that
$D_n^{\partial}\leq Y_{n,L_n}^{\partial}$.
Always, $D_n^{\partial}\leq\sum_{k=1}^{L_n}Y_{n,k}^{\partial}$.
Assumption (B3) and Minkowski's inequality yield
$$\big(E[(Y_{n,k}^{\partial})^4]\big)^{1/4}\leq C\tau_k/T_n,
\qquad E\left[\left(\sum_{k=1}^{L_n}Y_{n,k}^{\partial}\right)^4\right]\leq C.$$
Moreover, by (B2-$r$) and (\ref{LnTn-ineq2}),
$$P((E_n^{\partial})^c)
\leq\sum_{l=1}^2P(\mathcal N_{T_n}^{n,l}-\mathcal N_{T_n-\underline\tau}^{n,l}=0)
\leq C(\underline\tau/h_n)^{-r}\longrightarrow0.$$
H\"older's inequality and (\ref{LnTn-ineq}) therefore imply
$$E[D_n^{\partial}]
\leq C\bar\tau/T_n+C P((E_n^{\partial})^c)^{3/4}\longrightarrow0.$$
Consequently, $E[\|R_{n,11}\|]+E[\|R_{n,12}\|]\to0$.

Therefore, Lemma~\ref{Fn-lem} yields
\EQNN{\bar{E}[F_{n,11}]
&=\bar{E}\bigg[\SUMPQ \SUMLL \SUMMM \dot{\mca}_{l_1,l_2,m_1,m_2,p,q} \dot{F}_{l_1,l_2,m_1,m_2,p,q,n}\bigg]
+\bar E[R_{n,11}]\TOP0
}
as $n\to\infty$.
We similarly obtain $\bar{E}[F_{n,12}]\TOP 0$ and $\bar{E}[F_{n,22}]\TOP 0$,
and (B5) then yields $F_n\TOP\Gamma_1$.

\qed





\subsection{Proof of Proposition~\ref{renewal-exp-prop}}\label{Prop33-proof-subsection}

In this subsection, we prove Proposition~\ref{renewal-exp-prop}. 
To this end, we verify each condition (B1-q), (B2-q), (B3), (B4-$q$), (B5) and (B6). First, the following three lemmas hold.

\begin{lemma}\label{renewal-B1B2-lem}
Assume (R).
Then {\rm (B1-$q$)} and {\rm (B2-$q$)} are satisfied for any $q>0$.
\end{lemma}

\begin{proof}
For $0\leq t\leq (n-1)h_n$, we have
\[
\mcn_{t+h_n}^{n,l}-\mcn_t^{n,l}
=\bar{\mcn}_{h_n^{-1}t+1}^l-\bar{\mcn}_{h_n^{-1}t}^l.
\]
Since every inter-arrival time of $\bar{\mcn}^l$ is bounded from below by $\underline{\zeta}$,
the number of renewal epochs in any interval of length one is bounded by
$1+\underline{\zeta}^{-1}$. Therefore,
\[
\sup_{n\geq 1}\max_{l\in\{1,2\}}\sup_{0\leq t\leq (n-1)h_n}
E[(\mcn_{t+h_n}^{n,l}-\mcn_t^{n,l})^q]
\leq (1+\underline{\zeta}^{-1})^q<\infty,
\]
which proves (B1-$q$).

Next, since every inter-arrival time of $\bar{\mcn}^l$ is
bounded from above by $\bar{\zeta}$, a renewal epoch must occur within time
$\bar{\zeta}$ after any time. Hence
\[
P(\mcn_{t+uh_n}^{n,l}-\mcn_t^{n,l}=0)
=P(\bar{\mcn}_{h_n^{-1}t+u}^{l}-\bar{\mcn}_{h_n^{-1}t}^{l}=0)=0
\]
for $0\leq t\leq nh_n-uh_n$
whenever $u>\bar{\zeta}$, which proves (B2-$q$).
\end{proof}


\begin{lemma}\label{renewal-B3-prop}
Assume (R).
Then Assumption {\rm (B3)} is satisfied.
\end{lemma}

\begin{proof}

Recall that the sampling times are the renewal epochs and hence
$|I_i^l|=S_i^{n,l}-S_{i-1}^{n,l}=h_n\zeta_i^l$ for all $i,l$, so that
\[
r_n=\max_{i,l}|I_i^l|\leq \bar{\zeta}h_n
\quad\text{and}\quad
\underbar{r}_n=\min_{i,l}|I_i^l|\geq \underline{\zeta}h_n.
\]

Fix $\delta\in(0,1)$ and set $L_n=[n^{1-\delta}]$ and $b_n=[n/L_n]$.
Let $\ell_n=n-L_nb_n$ and define a partition
$(t_k^n)_{k=0}^{L_n}$ of $[0,T_n]$ by $t_0^n=0$ and
\[
t_k^n-t_{k-1}^n=h_n\big(b_n+1_{\{k\leq \ell_n\}}\big),\qquad k=1,\ldots,L_n,
\]
so that $t_{L_n}^n=nh_n=T_n$ and
\[
\underline{\tau}:=\min_{1\leq k\leq L_n}(t_k^n-t_{k-1}^n)=h_nb_n,\qquad
\bar{\tau}:=\max_{1\leq k\leq L_n}(t_k^n-t_{k-1}^n)=h_n(b_n+1).
\]
Hence $\sup_n \bar{\tau}/\underline{\tau}\leq 2$, $L_n\to\infty$, and, for any $\epsilon\in(0,\delta)$,
\[
n^{\epsilon-1}L_n\leq n^{\epsilon-\delta}\to 0.
\]
Moreover,
\[
L_n(nh_n)^{-1}E[r_n^4]^{1/4}
\leq L_n(nh_n)^{-1}\bar{\zeta}h_n
=\bar{\zeta}\,\frac{L_n}{n}\to 0.
\]

It remains to verify the moment bound in (\ref{B3-eq}).
Fix $n$ and $k$ and write $\tau_k=t_k^n-t_{k-1}^n$.
For each $l\in\{1,2\}$, the number of indices $i$ such that
$I_i^l\cap [t_{k-1}^n,t_k^n)\neq\emptyset$ is at most
$2+\tau_k/(\underline{\zeta}h_n)$ because every such interval has length at least
$\underline{\zeta}h_n$.
Therefore,
\[
\sum_{\substack{1\leq i\leq \M_n \\ I_i\cap [t_{k-1}^n,t_k^n)\neq\emptyset}}|I_i|^{-1}
\leq \sum_{l=1}^2\Big(2+\frac{\tau_k}{\underline{\zeta}h_n}\Big)\frac{1}{\underline{\zeta}h_n}
\leq \frac{C}{h_n^2}\,\tau_k
\]
with a constant $C>0$ independent of $n$ and $k$.
Consequently,
\[
\Big(\frac{h_n^2}{\tau_k}\sum_{\substack{1\leq i\leq \M_n \\ I_i\cap [t_{k-1}^n,t_k^n)\neq\emptyset}}|I_i|^{-1}\Big)^4
\leq C^4,
\]
and hence (\ref{B3-eq}) holds. This proves (B3).

\end{proof}

\noindent
{\bf Proof of Proposition~\ref{renewal-exp-prop}.}

By Theorem~1 and Corollary~1 of Bardet, Christen and
Fontbona~\cite{bar-etal15}, the residual-life transition laws satisfy
an exponential total variation bound, uniformly over initial
residual times in $(0,\bar\zeta]$.
Since the conditional law of the future renewal increments
depends on the past only through the residual times, the
coupling inequality and the independence of the two renewal
sequences yield
\[
\alpha_k^n\leq Ce^{-ck},\qquad k\geq1,
\]
for constants $C,c>0$ independent of $n$.
Consequently, (B4-$q$) holds for every $q>0$.

Then, by Lemmas~\ref{renewal-B1B2-lem} and~\ref{renewal-B3-prop}, we only need to show (B5) and (B6).
Moreover, by Lemma~\ref{Fn-lem}, to show (B5), we only need to show that there exist $\dot{F}_{l_1,l_2,m_1,m_2,p,q}$, $\check{F}_{l_1,l_2,m_1,m_2,p,q}$ and $\ddot{F}_{l_1,l_2,p}$ such that 
\EQ{\label{exp-F-conv} E[\DFN] \to \dot{F}_{l_1,l_2,m_1,m_2,p,q}, \quad E[\CFN]\to \check{F}_{l_1,l_2,m_1,m_2,p,q},
\quad E[\DDFN]\to \ddot{F}_{l_1,l_2,p}}
as $n\to\infty$.

The main idea for proving these expectation limits is to exploit the regenerative structure of the renewal sampling scheme.
For each local contribution $X'_k$ to $\DFN$, $\CFN$ or $\DDFN$, we first express $X'_k$ as $h_nT_n^{-1}$ times a function of the forward recurrence times after a deterministic grid point and of the subsequent renewal increments.
The latter increments are independent of the past and have the same distribution for all $k$, while the distribution of the forward recurrence times converges to the stationary excess-life distribution by the renewal theorem.
This gives the convergence of $E[X'_k]$ away from the boundary, and the boundary terms are negligible because the inter-arrival times are bounded.
Summing these local limits yields (B5); the same localization argument applied to the quantities in (B6) gives the required average conditions.

Let 
$s_{k,p,q}=k-1-2(p\vee q+1)\bar{\zeta}$, and let
$X'_k$ be one of $\CFK$, $\DDFK$, or $\DFK$.

Let 
\EQNN{\bar{\tau}_k^l=\min\{t>s_{k,p,q}; \Delta \bar{\mcn}_t^l>0\}, \quad
\tau_k^{n,l}=\min\{t>s_{k,p,q}h_n; \Delta \mcn_t^{n,l}>0\}.}
Let $\mcx_k=(\mcx_{1,k},\mcx_{2,k})=(\bar{\tau}_k^1-s_{k,p,q}, \bar{\tau}_k^2-s_{k,p,q})$
(forward recurrence time), then,
$$(\tau_k^{n,1}-s_{k,p,q}h_n, \tau_k^{n,2}-s_{k,p,q}h_n)=h_n\mcx_k.$$

\EQ{\label{renewal-exp-eq0} E[X'_k]=\int \int E[X'_k|\mcx_k=(s,t)]P\mcx_k^{-1}(dsdt).}


The first observation times after $s_{k,p,q}h_n$ are $s_{k,p,q}h_n+\mcx_{1,k}h_n$ and $s_{k,p,q}h_n+\mcx_{2,k}h_n$.
For $l\in \{1,2\}$, let $L_{l,k}=\bar{\mcn}^l_{\bar{\tau}_k^l}$, and 
$$U_{l,k}=\max\{\bar{\mcn}_t^l;t<k+2(p\vee q+1)\bar{\zeta}\}.$$

Let
$$\zeta^{(k)}=((\zeta_{i+L_{1,k}}^1)_{i=1}^\infty, (\zeta_{j+L_{2,k}}^2)_{j=1}^\infty).$$

For $\xi = ((\xi_i^1)_{i=1}^\infty, (\xi_j^2)_{j=1}^\infty)$ and $(s_1,s_2)\in \mbbr^2$, let $\varphi_0^l((s_1,s_2),\xi)=s_l$, and
$$\varphi_i^l((s_1,s_2),\xi)=s_l+\sum_{j=1}^i \xi_j^l \quad (i\geq 1).$$ 
Then, we can write
$$I_{i+L_{l,k}}^l=[s_{k,p,q}h_n+h_n\varphi_{i-1}^l, s_{k,p,q}h_n+h_n\varphi_i^l)(\mcx_k, \zeta^{(k)}), \quad |I_{i+L_{l,k}}^l|=h_n\xi_i^l.$$

We first consider the case that $X'_k = \ddot{F}_{1,1,p,n,k}$.
Let 
$$f_{i,j,i'}((s_1,s_2),\xi)=\frac{\big|[\varphi_{i-1}^1,\varphi_i^1)\cap [\varphi_{j-1}^2,\varphi_j^2)\big|((s_1,s_2),\xi)}{\sqrt{\xi_i^1}\sqrt{\xi_j^2}}
\times \frac{\big|[\varphi_{i'-1}^1,\varphi_{i'}^1)\cap [\varphi_{j-1}^2,\varphi_j^2)\big|((s_1,s_2),\xi)}{\sqrt{\xi_j^2}\sqrt{\xi_{i'}^1}},$$
Since $I_{i'}^1\cap [(k-1)h_n,kh_n)\neq \emptyset$, $I_{i_{2p'-2}}^1\cap I_{i_{2p'-1}}^2\neq \emptyset$ and
$I_{i_{2p'-1}}^2\cap I_{i_{2p'}}^1\neq \emptyset$ $(1\leq p'\leq p)$ imply that
$L_{2,k}<i_{2p'-1}<U_{2,k}$ $(1\leq p'\leq p)$ and $L_{1,k}<i_{2p'}<U_{1,k}$ $(0\leq p'\leq p)$.
Then, we have
\EQNN{\ddot{F}_{1,1,p,n,k}&= T_n^{-1} (\mfi_1^{(2)})^\top S_{1,1,p}\mfi_{1}^{(2,k)} \\
&=\frac{1}{T_n}\sum_{i_0,\cdots, i_{2p}}\bigg\{\prod_{p'=1}^p \frac{|I_{i_{2p'-2}}^1\cap I_{i_{2p'-1}}^2|}{\sqrt{|I_{i_{2p'-2}}|}\sqrt{|I_{i_{2p'-1}}|}}
\frac{|I_{i_{2p'-1}}^2\cap I_{i_{2p'}}^1|}{\sqrt{|I_{i_{2p'-1}}|}\sqrt{|I_{i_{2p'}}|}}\bigg\} \sqrt{|I_{i_0}^1|} \frac{|I_{i_{2p}}^1\cap [(k-1)h_n,kh_n)|}{\sqrt{|I_{i_{2p}}^1|}} \\
&=\frac{1}{T_n}\sum_{i'_0,\cdots, i'_{2p}=1}^\infty\bigg\{\prod_{p'=1}^p  f_{i'_{2p'-2},i'_{2p'-1},i'_{2p'}}(\mcx_k,\zeta^{(k)})\bigg\}
\frac{\sqrt{\zeta_{i'_0+L_{1,k}}^1}}{\sqrt{\zeta_{i'_{2p}+L_{1,k}}^1}} \\
&\qquad \times h_n\big|[\varphi_{i'_{2p}-1}^1, \varphi_{i'_{2p}}^1)(\mcx_k,\zeta^{(k)}) \cap [(2p+2)\bar{\zeta}-1, (2p+2)\bar{\zeta})\big| \\
&=\frac{h_n}{T_n}g(\mcx_k,\zeta^{(k)})
}
if 
\EQ{\label{k-cond-eq} (2p+2)\bar{\zeta}+1\leq k\leq n-(2p+2)\bar{\zeta},}
where $i'_{2p'-1}=i_{2p'-1}-L_{2,k}$, $i'_{2p'}=i_{2p'}-L_{1,k}$, and
\EQNN{g((s_1,s_2),\xi)&=\sum_{i_0,\cdots, i_{2p}=1}^\infty 
\bigg\{\prod_{p'=1}^p f_{i_{2p'-2},i_{2p'-1},i_{2p'}}((s_1,s_2),\xi)\bigg\}
\frac{\sqrt{\xi_{i_0}^1}}{\sqrt{\xi_{i_{2p}}^1}} \\
&\qquad \times \big|[\varphi_{i_{2p}-1}^1, \varphi_{i_{2p}}^1)((s_1,s_2),\xi) \cap [(2p+2)\bar{\zeta}-1, (2p+2)\bar{\zeta})\big|.
}
\begin{note}
{\colorr $g$ seems not to be able to be defined when $\xi$ is equal to zero for sufficiently large index.}
\end{note}

For fixed p and q, only observation intervals reachable from the unit cell through at most 2 max(p,q)+1 overlaps can contribute to g. Since the renewal intervals are bounded above and bounded away from zero, all such intervals lie in a fixed bounded region and their indices are uniformly bounded. Hence, g depends on only finitely many coordinates of the renewal increment sequences, with a bound independent of the phase variables and the increment sequences.

Note that $g$ does not depend on $n$ nor $k$, and is continuous on $(s_1,s_2)$.
Moreover, by a similar estimate to Lemma~\ref{three-F-lem}, $g(\mcx_k, \zeta^{(k)})$ is bounded uniformly on $k$ and $n$.

In the other cases of $X'_k$, we may write
\EQ{\label{xk-g-eq} X'_k = \frac{h_n}{T_n} g(\mcx_k, \zeta^{(k)})}
for a function $g$ that possesses similar properties as above.

Fix $k$ satisfying (\ref{k-cond-eq}).
Then, we have
\EQNN{E[X'_k \mid \mcx_k=(s,t)]
&=E[E[X'_k \mid L_{1,k},L_{2,k},\mcx_k]\mid \mcx_k=(s,t)] \\
&= h_nT_n^{-1}E[E[g((s,t),((\zeta_{L_{1,k}+1}^1,\cdots),(\zeta_{L_{2,k}+1}^2,\cdots))) \mid L_{1,k},L_{2,k},\mcx_k] \mid \mcx_k=(s,t)]
}

For any $A\in \sigma(L_{1,k},L_{2,k},\mcx_k)$ and $m_1, m_2\in \mbbz_+$, 
$g((s,t),((\zeta_{m_1+1}^1,\cdots),(\zeta_{m_2+1}^2,\cdots)))$ is\\
$\sigma(\zeta_{m_1+1},\cdots,\zeta_{m_2+1}^2,\cdots)$-measurable while 
$1_{A\cap \{(L_{1,k},L_{2,k})=(m_1,m_2)\}}$ is $\sigma(\zeta_1^1,\cdots, \zeta_{m_1}^1, \zeta_1^2,\cdots, \zeta_{m_2}^2)$-\\
measurable.
Then, these two random variables are independent, and hence we have
\EQNN{&E[g((s,t),((\zeta_{L_{1,k}+1}^1,\cdots),(\zeta_{L_{2,k}+1}^2,\cdots)))1_A] \\
&=\sum_{m_1,m_2=0}^\infty E[g((s,t),((\zeta_{m_1+1}^1,\cdots),(\zeta_{m_2+1}^2,\cdots)))1_{A\cap \{(L_{1,k},L_{2,k})=(m_1,m_2)\}}] \\
&=\sum_{m_1,m_2=0}^\infty E[h(s,t)1_{A\cap \{(L_{1,k},L_{2,k})=(m_1,m_2)\}}] \\
&=E[h(s,t)1_A],
}
where $h(s,t)=E[g((s,t),((\zeta_{1}^1,\cdots),(\zeta_{1}^2,\cdots)))]$.
Therefore, we obtain
$$E[g((s,t),((\zeta_{L_{1,k}+1}^1,\cdots),(\zeta_{L_{2,k}+1}^2,\cdots))) \mid L_{1,k},L_{2,k},\mcx_k] = h(s,t) \quad {\rm a.s.},$$
which implies
\EQ{\label{renewal-exp-eq1} E[X'_k|\mcx_k=(s,t)]=h_nT_n^{-1}h(s,t) \quad P\mcx_k^{-1}\mathchar`-{\rm a.e.} (s,t).}

By the properties of $g$, we obtain that $h$ is continuous and bounded with respect to  $(s,t)$ by the dominated convergence theorem. 
Moreover, Theorem 6.10 in Medhi ``Stochastic Processes'' P187 yields that 
\EQ{\label{renewal-exp-eq2} P\mcx_k^{-1}\Rightarrow \frac{1-F_1(s)}{E[\zeta_1^1]}\frac{1-F_2(t)}{E[\zeta_1^2]}dsdt =:f(s,t)dsdt.}
where $F_l$ is the distribution function of $\zeta_1^l$.

(\ref{renewal-exp-eq0}), (\ref{renewal-exp-eq1}), (\ref{renewal-exp-eq2}) yield that there exists $N$ such that $k\geq N$ implies
\EQ{\label{A6EXk-eq} \bigg|E[X'_k]-\frac{h_n}{T_n}\int \int h(s,t)f(s,t)dsdt\bigg|<\frac{h_n\epsilon}{T_n}.}

We obtain $\tau_{k,p}\leq C_ph_n$ by $\xi_k^l\leq \bar{\zeta}$. Then, Lemma~\ref{three-F-lem} yields $E[X'_k]\leq C_pn^{-1}$.
Therefore, we can ignore finite numbers of $k$ not satisfying (\ref{k-cond-eq}), and we have
$$E\bigg[\sum_{k=1}^nX'_k\bigg]\to \int \int h(s,t)f(s,t)dsdt$$
as $n\to\infty$, which implies (\ref{exp-F-conv}),
and hence (B5) holds.

To show (B6), we consider the mixed term; the other terms are treated
similarly. Fix $p\in\mathbb Z_+$ and a deterministic partition
$(s_k)_{k\geq0}\in\mathfrak S$. For an interval $A\subset[0,T_n]$, set
\[
 Q_n(A)=\sum_i\frac{|I_i^1\cap A|}{\sqrt{|I_i^1|}}
             [(GG^\top)^pG\mathfrak I_2]_i.
\]
For $C_{n,j}=[(j-1)h_n,jh_n)$, the preceding reconstruction and renewal
argument give, outside a fixed number of boundary cells,
\[
 Q_n(C_{n,j})=h_ng_j,\qquad |g_j|\leq C_p,\qquad
 E[g_j]\longrightarrow f_p^{1,2},
\]
where the law of $g_j$ does not depend on $n$.
Let $B_k=[s_{k-1},s_k)$ and let $J_{n,k}$ index the interior cells fully
contained in $B_k$. Additivity and the local bound yield
\[
 Q_n(B_k)=h_n\sum_{j\in J_{n,k}}g_j+R_{n,k},\qquad
 \sup_{k\leq q_n}|R_{n,k}|=O(h_n).
\]
The remainder contains the partial endpoint cells and the excluded
boundary cells. Split each sum at a sufficiently large fixed index $N$.
The initial part of its expectation error is $O(Nh_n)$; the remaining
part is bounded by $\sup_{j\geq N}|E[g_j]-f_p^{1,2}|$, uniformly in $k$
since $|B_k|\leq1$. Letting $n\to\infty$ and then $N\to\infty$ gives
uniform convergence to $f_p^{1,2}|B_k|$.

Replacing $|I_i^1\cap B_k|/|I_i^1|$ by
$1_{\{\sup I_i^1\in B_k\}}$ affects at most two intervals.
Since $(GG^\top)^pG\mathfrak I_2\leq\mathfrak I_1$ componentwise,
the resulting error is at most $2\bar\zeta h_n$. Consequently,
\[
 \max_{1\leq k\leq q_n}
 \left|E[\mathfrak I_1^\top\mathcal E_{(k)}^1
              (GG^\top)^pG\mathfrak I_2]
       -f_p^{1,2}(s_k-s_{k-1})\right|\longrightarrow0.
\]
The other limits follow similarly using the corresponding counting,
diagonal-trace and bilinear local quantities; the counting constants are
$a_0^l=1/E[\zeta_1^l]$ by the renewal theorem.
The residual-block boundedness follows from the interval-length bounds.
\begin{note}
{\colorr Boundedness of $g$ is derived by $\|(GG^\top)^p\|\leq 1$ and $\zeta\geq \underline{\zeta}$.}
\end{note}

Finally, since $a_0^l=1/E[\zeta_1^l]$ and
$0<\underline{\zeta}\leq\zeta_1^l\leq\bar{\zeta}<\infty$
almost surely, we have
\[
0<\bar{\zeta}^{-1}\leq a_0^l\leq\underline{\zeta}^{-1}<\infty,
\qquad l=1,2.
\]
This shows that $a_0^1$ and $a_0^2$ are positive.

\qed


\begin{small}
\bibliographystyle{abbrv}
\bibliography{reference}
\end{small}

\renewcommand{\thesection}{\Alph{section}}
\setcounter{section}{0}
\section{Some supplementary results}

\begin{lemma}\label{alpha-k-lem}
Let $u\geq 0$, $\beta\geq 1$, $m$ be a positive integer, and $(a_k)_{k=0}^\infty$ be a monotone non-increasing sequence.
Let $K=\inf_{k\geq 1}([k\beta]-[(k-1)\beta])$. Then, $K\geq 1$ and 
$$\sum_{k=1}^\infty1_{\{a_{([k\beta]-m)\vee 0}>u\}}
\leq K^{-1}\bigg(\sum_{k=0}^\infty 1_{\{a_k>u\}}+m-1\bigg).$$
\end{lemma}

\begin{proof}
Since $(a_k)_{k=0}^\infty$ is monotone,
We have $1_{\{a_{([k\beta]-m)\vee 0}>u\}}\leq 1_{\{a_{k'\vee 0}>u\}}$
for a positive integer $k$ and $[(k-1)\beta]-m+1\leq k'\leq [k\beta]-m$,.
Therefore, we obtain
\EQNN{\sum_{k=1}^\infty1_{\{a_{([k\beta]-m)\vee 0}>u\}} 
&\leq \sum_{k=1}^\infty \frac{1}{[k\beta]-[(k-1)\beta]}\sum_{k'=[(k-1)\beta]-m+1}^{[k\beta]-m}1_{\{a_{k'\vee 0}>u\}} \\
&\leq K^{-1}\sum_{k'=-m+1}^\infty 1_{\{a_{k'\vee 0}>u\}} \\
&\leq K^{-1}\bigg(\sum_{k'=0}^\infty 1_{\{a_{k'}>u\}}+m-1\bigg).
}

\end{proof}

\noindent
{\bf Proof of Lemma~\ref{GGTp-lem}.}

\EQNN{[(GG^\top)^pG]_{ii'}\frac{1}{\sqrt{|I_{i'}|}}
=\sum_{\substack{i_1,\cdots, i_p \\ j_1,\cdots, j_p}}
\bigg(\prod_{l=1}^p\frac{|I_{i_{l-1}}\cap I_{j_l}|}{\sqrt{|I_{i_{l-1}}|} \sqrt{|I_{j_l}|}}\frac{|I_{j_l}\cap I_{i_l}|}{\sqrt{|I_{j_l}|}\sqrt{|I_{i_l}|}}\bigg)
\cdot \frac{|I_{i_p}\cap I_{i'}|}{\sqrt{|I_{i_p}|}\sqrt{|I_{i'}|}}\cdot \frac{1}{\sqrt{|I_{i'}|}}.
}

Since 
$$\sum_{j_l}\frac{|I_{i_{l-1}}\cap I_{j_l}|}{\sqrt{|I_{i_{l-1}}|}\sqrt{|I_{j_l}|}} \frac{|I_{j_l}\cap I_{i_l}|}{\sqrt{|I_{j_l}|}\sqrt{|I_{i_l}|}}
\leq \frac{1}{\sqrt{|I_{i_{l-1}}|}}\sum_{j_l}\frac{|I_{j_l}\cap I_{i_l}|}{\sqrt{|I_{i_l}|}} = \frac{\sqrt{|I_{i_l}|}}{\sqrt{|I_{i_{l-1}}|}},$$
$$\sum_{j_l}\sqrt{|I_{j_l}|} \frac{|I_{j_l}\cap I_{i_l}|}{\sqrt{|I_{j_l}|}\sqrt{|I_{i_l}|}}
\leq \sqrt{|I_{i_l}|},
\quad \sum_{i_l}\sqrt{|I_{i_l}|} \frac{|I_{i_l}\cap I_{j_{l+1}}|}{\sqrt{|I_{i_{l}}|}\sqrt{|I_{j_{l+1}}|}}
\leq  \sqrt{|I_{j_{l+1}}|},$$
we have
$$[(GG^\top)^pG]_{ii'}\frac{1}{\sqrt{|I_{i'}|}}
\leq \frac{1}{\sqrt{|I_i|}}\cdot \sqrt{|I_{i'}|}\cdot \frac{1}{\sqrt{|I_{i'}|}}=\frac{1}{\sqrt{|I_i|}}.$$

We similarly obtain the inequality for $[(G^\top G)^pG^\top]_{i'i}$.

\qed

\noindent
{\bf Proof of Lemma~\ref{three-F-lem}.}

(\ref{mca-eq}) is obtained by (A1), (A2) and 
$|\PB^l a_{l_1,l_2,p}|\leq C(p+1) \bar{\rho}^{(2p-l)\vee 0}$
for $l\in \{0,1\}$ and $p\in \mbbz_+$.
Let $\mci_k^{l,p}$ be the set of all index $1\leq i\leq \M_{l,n}$
such that for some $i_1,\cdots, i_p\in \{1,\cdots, \M_n\}$,
$I_{i_{j-1}}\cap I_{i_j}\neq \emptyset$ $(2\leq j\leq p)$, $I_i^l\cap I_{i_1}\neq \emptyset$
and $I_{i_p}\cap [t_{k-1},t_k)\neq \emptyset$.
Let $\mfi_l^{(2,q,k)}=(|I_i^l|1_{\{i\in \mci_k^{l,q}\}})_{i=1}^{\M_{l,n}}$.

Then, (\ref{G-ineq}) yields
$$|\ddot{F}_{l_1,l_2,q,n,k}| \leq \frac{1}{T_n}|D_{l_1}^{-1/2}\mfi_{l_1}^{(2,2q+1,k)}|\sqrt{\sum_{i\in \mci_k^{l_2}}|I_i|}
\leq T_n^{-1}\tau_{k,q},$$
and hence we have
$$E[\ddot{F}_{l_1,l_2,q,n,k}^4]\leq T_n^{-4}E[\tau_{k,q}^4]\leq CT_n^{-4}(\tau_k^4+(q^4+1)E[r_n^4]).$$

Moreover, Lemma~\ref{GGTp-lem} and (\ref{G-ineq}) yield
\EQNN{|\dot{F}_{1,2,1,2,p,q,n,k}|
&\leq \frac{h_n}{n}\sum_{i,i'}\frac{|I_i\cap I_{i'}\cap \TKI|}{|I_i||I_{i'}|}[(GG^\top )^pG]_{ii'}[(GG^\top)^qG]_{ii'} \\
&\leq \frac{h_n}{n}\sum_{i,i'}\frac{|I_i\cap I_{i'}\cap \TKI|}{|I_i|^2} \\
&\leq \frac{h_n}{n}\sum_{i\in \mci_k^1}|I_i^1|^{-1},
}
$$|\dot{F}_{1,1,1,1,p,q,n,k}|
\leq \frac{h_n}{n}\sum_{i,i'}|I_i\cap I_{i'}\cap \TKI|\frac{[(GG^\top)^p]_{ii}[(GG^\top)^q]_{ii}}{|I_i|^2}
\leq \frac{h_n}{n}\sum_i\frac{|I_i\cap \TKI|}{|I_i|^2}
\leq \frac{h_n}{n}\sum_{i\in \mci_k^1}|I_i^1|^{-1},$$
\EQNN{|\dot{F}_{1,1,1,2,p,q,n,k}|
&\leq \frac{h_n}{n}\sum_{i,i'}|I_i\cap I_{i'}\cap \TKI| \frac{[(GG^\top )^p]_{ii}[(GG^\top)^qG]_{ii'}}{|I_i|^{3/2} \sqrt{|I_{i'}|}} \\
&\leq \frac{h_n}{n}\sum_{i\in \mci_k^1,i'}\frac{|I_i\cap I_{i'}|}{\sqrt{|I_i|}\sqrt{|I_{i'}|}} \frac{[(GG^\top)^qG]_{ii'}}{|I_i|} \\
&= \frac{h_n}{n}\sum_{i\in \mci_k^1} \frac{[(GG^\top)^{q+1}]_{ii}}{|I_i|} \\
&\leq \frac{h_n}{n}\sum_{i\in \mci_k^1}|I_i^1|^{-1}.
}
Furthermore, the Cauchy-Schwarz inequality yields
\EQNN{|\dot{F}_{1,1,2,2,p,q,n,k}|
&\leq \frac{h_n}{n}\sum_{i,i'}\frac{|I_i\cap I_{i'}\cap \TKI|}{|I_i||I_{i'}|}[(GG^\top)^p]_{ii}[(G^\top G)^q]_{i'i'} \\
&\leq \frac{h_n}{n}\sum_{i,i'}\frac{\sqrt{|I_i\cap I_{i'}\cap \TKI|}}{|I_i|}\frac{\sqrt{|I_i\cap I_{i'}\cap \TKI|}}{|I_{i'}|} \\
&\leq \frac{h_n}{n}\sqrt{\sum_{i\in \mci_k^1,i'}\frac{|I_i\cap I_{i'}|}{|I_i|^2}}\sqrt{\sum_{i,i'-\M_{1,n}\in \mci_k^2}\frac{|I_i\cap I_{i'}|}{|I_{i'}|^2}} \\
&\leq \frac{h_n}{n}\sqrt{\sum_{i\in \mci_k^1}|I_i|^{-1}} \sqrt{\sum_{i'-\M_{1,n}\in\mci_k^2}|I_{i'}|^{-1}} \\
&\leq \frac{h_n}{n}\sum_{i\in \mci_k}|I_i|^{-1}.
}
Similar estimates yield
$$|\DFK| \leq \frac{h_n}{n} \sum_{i\in \mci_k}|I_i|^{-1}$$
for any $l_1,l_2,m_1,m_2,p,q,n,k$.
Moreover, (B3) implies that
$$E[\DFK^4]\leq Ch_n^4/n^4  \cdot \tau_k^4h_n^{-8}=CT_n^{-4}\tau_k^4.$$

\EQNN{|\check{F}_{1,1,1,m_2,p,q,n,k}|
&\leq \frac{1}{n}\sum_{i\in \mci_k^1, i'} \frac{|I_i\cap I_{i'}|}{|I_i|^{3/2}}[(GG^\top)^p]_{ii}[D_1^{1/2}S_{1,m_2,q}\mfi_{m_2}^{(2)}]_i \\
&\leq \frac{1}{n}\sum_{i\in \mci_k^1}|I_i|^{-1/2}[D_1^{1/2}S_{1,m_2,q}\mfi_{m_2}^{(2,2q+1,k)}]_i \\
&\leq \frac{1}{n}\sqrt{\sum_{i\in \mci_k^1}|I_i|^{-1}}\lVert D_1^{1/2}S_{1,m_2,q} D_{m_2}^{1/2}\rVert |D_{m_2}^{-1/2}\mfi_{m_2}^{(2,2q+1,k)}| \\
&\leq \frac{1}{n}\sqrt{\tau_{k,q}}\sqrt{\sum_{i\in \mci_k^1}|I_i|^{-1}},
}
\EQNN{|\check{F}_{1,1,2,m_2,p,q,n,k}|
&\leq \frac{1}{n}\sum_{i\in \mci_k^1, i'}\frac{|I_i\cap I_{i'}|}{|I_i|\sqrt{|I_{i'}|}}[(GG^\top)^p]_{ii}[D_2^{1/2}S_{2,m_2,q}\mfi_{m_2}^{(2)}]_{i'} \\
&\leq \frac{1}{n}\sum_{i\in \mci_k^1, i'}|I_i|^{-1/2}[G]_{ii'}[D_2^{1/2}S_{2,m_2,q}\mfi_{m_2}^{(2)}]_{i'} \\
&\leq \frac{1}{n}\sqrt{\sum_{i\in \mci_k^1} |I_i|^{-1}}\lVert G\rVert \lVert D_2^{1/2}S_{2,m_2,q}D_{m_2}^{1/2} \rVert |D_{m_2}^{-1/2}\mfi_{m_2}^{(2,2q+2,k)}| \\
&\leq \frac{1}{n}\sqrt{\tau_{k,q}}\sqrt{\sum_{i\in \mci_k^1} |I_i|^{-1}},
}
\EQNN{|\check{F}_{1,2,1,m_2,p,q,n,k}|
&\leq \frac{1}{n}\sum_{i\in \mci_k^1, i'}\frac{|I_i\cap I_{i'}|}{|I_i|\sqrt{|I_{i'}|}} [(GG^\top)^pG]_{ii'}[D_1^{1/2}S_{1,m_2,q}\mfi_{m_2}^{(2)}]_i \\
&=\frac{1}{n}\sum_{i\in \mci_k^1}|I_i|^{-1/2}[(GG^\top)^{p+1}]_{ii} [D_1^{1/2}S_{1,m_2,q}\mfi_{m_2}^{(2)}]_i \\
&\leq \frac{1}{n}\sqrt{\sum_{i\in\mci_k^1}|I_i|^{-1}}\lVert D_1^{1/2}S_{1,m_2,q} D_{m_2}^{1/2}\rVert |D_{m_2}^{-1/2}\mfi_{m_2}^{(2,2q+1,k)}| \\
&\leq \frac{1}{n}\sqrt{\tau_{k,q}}\sqrt{\sum_{i\in\mci_k^1}|I_i|^{-1}}.
}
Other terms can be estimated in similar ways.

Finally, (B3) and the Cauchy-Schwarz inequality imply that
\EQNN{E[\CFK^4] &\leq n^{-4} E[\tau_{k,q}^4]^{1/2} E\bigg[\bigg(\sum_{i\in \mci_k}|I_i|^{-1}\bigg)^4\bigg]^{1/2} \\
&\leq C n^{-4} (\tau_k^2 + (q^2+1)E[r_n^4]^{1/2}) \cdot \tau_k^2 h_n^{-4} \\
&=CT_n^{-4}(\tau_k^4+(q^4+1)E[r_n^4]).
}

\qed




\end{document}